\documentclass[11pt]{article}
\usepackage{amssymb,latexsym}
\usepackage{amsmath}
\newcommand{\df}{\displaystyle\frac}

\newcommand{\hsp}{\hspace*{\parindent}}

\newcommand{\sig}{\sigma}

\newcommand{\sE}{{\cal E}}

\newcommand{\ra}{\rightarrow}
\newcommand{\In}{\infty}
\newcommand{\RR}{{\Bbb R}}
\newcommand{\ZZ}{{\Bbb Z}}
\newcommand{\NN}{{\Bbb N}}
\newcommand{\EE}{{\Bbb E}}
\newcommand{\QQ}{{\Bbb Q}}

\newcommand{\sR}{{\cal R}}

\newcommand{\sA}{{\cal A}}

\newcommand{\sS}{{\cal S}}
\newcommand{\sT}{{\cal T}}

\newcommand{\ba}{{\bf a}}

\newcommand{\by}{{\bf y}}

\input amssym.def
\input amssym.tex

\newcommand{\beql}[1]{\begin{equation}\label{#1}}
\newcommand{\eeq}{\end{equation}}

\catcode`\@=11
\renewcommand{\section}{
        \setcounter{equation}{0}
        \@startsection {section}{1}{\z@}{-3.5ex plus -1ex minus
        -.2ex}{2.3ex plus .2ex}{\large\bf}
        }
\catcode`@=12
\makeatletter
\def\eqalignno#1{\displ@y \tabskip\@centering
  \halign to\displaywidth{\hfil$\@lign\displaystyle{##}$\tabskip\z@skip
    &$\@lign\displaystyle{{}##}$\hfil\tabskip\@centering
    &\llap{$\@lign##$}\tabskip\z@skip\crcr
    #1\crcr}}
\makeatletter
\def\eqalignno#1{\displ@y \ta {\bf s} kip\@centering
  \halign to\displaywidth{\hfil$\@lign\displaystyle{##}$\ta {\bf s} kip\z@skip
  \ldots \@lign\displaystyle{{}##}$\hfil\ta {\bf s} kip\@centering
  \ldots llap{$\@lign##$}\ta {\bf s} kip\z@skip\crcr
    #1\crcr}}
\makeatother

\makeatletter
\def\@sect#1#2#3#4#5#6[#7]#8{\ifnum #2>\c@secnumdepth
     \def\@svsec{}\else 
     \refstepcounter{#1}\edef\@svsec{\csname the#1\endcsname.\hskip .75em }\fi
     \@tempskipa #5\relax
      \ifdim \@tempskipa>\z@ 
        \begingroup #6\relax
          \@hangfrom{\hskip #3\relax\@svsec}{\interlinepenalty \@M #8\par}%
        \endgroup
       \csname #1mark\endcsname{#7}\addcontentsline
         {toc}{#1}{\ifnum #2>\c@secnumdepth \else
                      \protect\numberline{\csname the#1\endcsname}\fi
                    #7}\else
        \def\@svsechd{#6\hskip #3\@svsec #8\csname #1mark\endcsname
                      {#7}\addcontentsline
                           {toc}{#1}{\ifnum #2>\c@secnumdepth \else
                             \protect\numberline{\csname the#1\endcsname}\fi
                       #7}}\fi
     \@xsect{#5}}
\def\@begintheorem#1#2{\it \trivlist \item[\hskip \labelsep{\bf #1\ #2.}]}
\makeatother

\begin{document}
\begin{center}
{\Large {\bf The $3x + 1$ Problem: An Annotated Bibliography, II (2000-2009)}} \\
\vspace{\baselineskip}
{\large {\em J. C. Lagarias}} \\
\vspace*{.2\baselineskip}
Department of Mathematics \\
University of Michigan \\
Ann Arbor, MI 48109--1109\\
{\tt lagarias@umich.edu}\\
\vspace*{1\baselineskip}
(Jan. 10, 2012) \\
\vspace*{2\baselineskip}
\end{center}

\noindent{\bf ABSTRACT.} 
The $3x+1$ problem concerns iteration of the map $T: \ZZ \rightarrow \ZZ$
given by
$$
T(x) = \left\{
\begin{array}{cl}
\df{3x+1}{2} & \mbox{if} ~~x \equiv 1~~ (\bmod~2) ~. \\
~~~ \\
\df{x}{2} & \mbox{if}  ~~x \equiv 0~~ (\bmod~2)~.
\end{array}
\right.
$$
The $3x+1$ Conjecture asserts that each $m \geq 1$ has some iterate
$T^{(k)} (m) = 1$.
This is the second installment of
an annotated bibliography of work done on the $3x+1$ problem
and related problems, mainly covering the period 2000 through 2009,
with some related later  papers (which were preprints by 2009).
At present the  $3x+1$ Conjecture remains unsolved. 

%
%
%
\setlength{\baselineskip}{1.0\baselineskip}
\section{Introduction}
\hsp

The $3x+1$ problem is most simply stated in terms
of the {\em Collatz function} $C(x)$ defined on integers as
``multiply by three and add one'' for odd integers
and ``divide by two'' for even integers. That is,
$$
C(x) = 
\left\{
\begin{array}{cl}
3x+1 & \mbox{if}~ x \equiv 1~~ (\bmod ~2 ) ~, \\
~~~ \\
\df{x}{2} & \mbox{if} ~~x \equiv 0~~ (\bmod~2) ~,
\end{array}
\right.
$$
The 
{\em $3x+1$ problem} (or {\em Collatz problem}) 
is to prove that starting from any positive integer,
some iterate of this function takes the value $1$.
The problem  other names: it
has also been called  Kakutani's problem, the Syracuse problem,
 and Ulam's problem.

Much work on the problem is stated in terms 
of the  {\em $3x+1$  function} 
$$
T(x) = \left\{
\begin{array}{cl}
\df{3x+1}{2}  & \mbox{if} ~~x \equiv 1~~ (\bmod ~ 2) \\
~~~ \\
\df{x}{2} & \mbox{if}~~ x \equiv 0~~ (\bmod ~2 )~.
\end{array}
\right.
$$
The  $3x  +  1$ {\em Conjecture} 
states that every $m \geq  1$ has some iterate
$T^{(k)} (m) = 1$.

The $3x+1$ Conjecture has  now been 
verified  up to  $17 \times 2^{58} >  4.899 \times 10^{18}$ 
(as of Feb. 21, 2008) by an ongoing computation 
run by  T. Oliveira e Silva (2004+).  An independent
computation of  Roosendaal(2004+) verifies it to $612 \times 2^{50}> 6.89  \times 10^{17}$.

At present the $3x+1$ conjecture  remains unsolved. 
The  proofs claimed in  Yamada (1981), Cadogan (2006) and Bruckman (2008)  are
incomplete. 

Surveys on results on the $3x+1$ problem can be found in
Lagarias (1985), M\"{u}ller (1991), and
the first chapter of Wirsching (1998a), described in
the first installment of the annotated bibliography, Lagarias (2003+).
A more recent survey appears in  Chamberland (2003).

%
%
%
\setlength{\baselineskip}{1.0\baselineskip}
\section{Terminology}
\hsp

We  use the following definitions.
The {\em trajectory} or {\em foward orbit} of an integer $m$ is the set
$$
O^{+}(m) :=\{ m,~ T(m)~, ~~ T^{(2)} (m) , \ldots \}~.
$$
The
{\em stopping time} $\sigma(m)$
of $m$ is the least $k$ such that $T^{(k)} (m) < m$,
and is $\infty$ if no such $k$ exists.
The
{\em total stopping time }
$\sig_\infty (m)$ is the least $k$ such that
$m$ iterates to $1$ under $k$ applications of the function  $T$ i.e.
$$ \sig_\infty (m) : = \inf~ \{k~:~ T^{(k)} (m) = 1\}. $$  
The {\em scaled total stopping time} or {\em gamma value} $\gamma (m)$ is given by
$$ \gamma (m) := \frac {\sig_\infty (m)} {log~ m}.$$
The
{\em height} $h(m)$ is 
the least $k$ for which the Collatz function $C(x)$ 
has $C^{(k)} (m) =  1$. It is also given by
$$h(m) := \sig_\infty (m) + d(m),$$ 
where $d(m)$ counts the
number of iterates $T^{(k)} (m) \equiv 1$ (mod~2) for
$0 \leq k < \sigma_\infty (m)$.
Finally, the  function $\pi_a (x)$ counts the number of $n$ with
$|n| \leq x$ that contain $a$ in their  forward orbit under $T$.

%
%
%

\section{Bibliography }
\hsp
This bibliography  covers research 
articles, survey articles and PhD theses on the $3x+1$ problem and
related problems from 2000  to the present. 
The first installment of the annotated bibliography is 
Lagarias(2003+), which covers the period 1963--1999.
Articles in Chinese  have the authors surname
listed first.

\begin{enumerate}

\item
Ethan Akin (2004), {\em Why is the $3x+1$ Problem Hard?,}
In: {\em Chapel Hill Ergodic Theory Workshops}  (I. Assani, Ed.), 
Contemp. Math. vol 356, Amer. Math. Soc. 2004, pp. 1--20. 
(MR 2005f:37031).\\
\newline
\hspace*{.25in}
This paper analyzes the $3x+1$ problem by viewing the 
map $T$ as acting on the
domain $\ZZ_2$ of 2-adic integers.
The map $T$ is topologically conjugate over $\ZZ_2$ 
to the 2-adic shift map
$$
S(x) = \left\{
\begin{array}{cl}
\df{x-1}{2} & \mbox{if}~~x \equiv 1~~ (\bmod ~2 ) ~, \\
~~~ \\
\df{x}{2} & \mbox{if}~~x \equiv 0~~(\bmod~2 ) ~,
\end{array}
\right.
$$
by a conjugacy map $Q_3 : \ZZ_2  \rightarrow  \ZZ_2$, 
i.e. $Q_3  \circ T  = S \circ Q_3$. (The map $Q_3$ equals
the map denoted $Q_{\infty}$ in Lagarias (1985), and is
the inverse of the map $\Phi$ in Bernstein (1994).)
The $3x+1$ Conjecture can be reformulated in terms of the behavior of 
$Q_{3}$ acting on integers, 
namely that $Q_3$ maps $\ZZ^+ $ into $\frac{1}{3} \ZZ$.
Consider  more generally for any odd rational $a$ the map $T_a(x)$
which sends $x \mapsto \frac{ax+1}{2}$ or $\frac{x}{2}$,
according as $x$ is an odd or even $2$-adic integer.
The author observes there  is an 
associated conjugacy map $Q_a:\ZZ_2  \rightarrow  \ZZ_2$ 
with the same property as above,  
and formulates the  Rationality Conjecture for $Q_a$,
which asserts that $Q_a$ maps the rationals
with odd denominators to rationals. He shows that
the Rationality conjecture is true for $a = \pm 1$ and
is false for any odd rational $a$ that is not an integer.
For the  remaining cases of odd integer $a$, where the
Rationality Conjecture remains unsolved, he presents
 a heuristic argument suggesting that 
it should be  true for $a = \pm 3$ 
and false for all odd integers $|a| \ge 5$.

\item
Jo\~{a}o F. Alves, M\'{a}rio M. Graca, M. E. Sousa Dias, and Jos\'{e} Sousa Ramos (2005), 
{\em A linear algebra approach to the conjecture of Collatz,}
Lin. Alg. Appl. {\bf 394} (2005), 277--289. (MR2100588) \\
\newline 
\hspace*{.25in}
This paper studies the conjecture that the only periodic orbit
of the Collatz map on the positive integers goes through $n=1$.
They form an $n \times n$ zero-one matrix $A_n$ whose entries are
$$
A_{i, j}= 1 ~~~\mbox{if}~~~T(i)=j, ~~1 \le i,  j \le n.
$$
where $T(n)$ is the $3x+1$ function, and $A_{i, j}=0$ otherwise.
The assertion that $\{1, 2\}$ is the only periodic  orbit  of $T$ on
the positive integers is shown to be equivalent to
$ \det(I - xA_n) = 1-x^2$ for all $n \ge 1$. They prove that
$\det(I-xA_n) = \det(I- xA_{n-1})$ for all $n \ne 8~(\bmod~18)$.
They deduce that if there is another periodic orbit on
the positive integers then there exists $m \equiv 8 ~(\bmod~  18)$
such that $n= \frac{m}{2}$ is in a periodic orbit. Various further
conditions are deduced in the case $n \equiv 8 (\bmod~18)$,
e.g. $\det(I-xA_n) = \det(I- xA_{n-1})$ if $n \equiv 8~(\bmod~54)$.

\item
Tewodrus Amdeberhan, Dante Manna and Victor H. Moll (2008),
{\em The $2$-adic valuation of a sequence arising from a rational integral,}
J. Combinatorial Theory, Series A, {\bf 115} (2008), no. 8, 1474--1486.
(MR 
2009k:11194) \\
\newline 
\hspace*{.25in}
This paper studies certain integer sequences $\{ A_{k, m}: k \ge 0\}$
arising from evaluation of the integral
$$
N_{0, 4}(a; m) = \int_{0}^{\infty} \frac{dx}{(x^4 + 4ax^2 +1)^{m+1}}
$$
expanded in Taylor series in the parameter $a$, as
$$
N_{0,4}(a;m) = \frac{\pi}{\sqrt{2} m! (4(2a+1))^{\frac{m+1}{2}}} 
\sum_{k=0}^{\infty} A_{k,m} \frac{ a^k}{k!}.
$$
These sequences are given by
$$
A_{k, m} = \frac{k!m!}{2^{m-k}}~ \sum_{j=k}^m 2^j 
\left({{2m-2j}\atop{m-j}}\right)
\left({{m+j}\atop{m}}\right)\left({{j}\atop{k}}\right).
$$

In section 6 a relation is shown between divisibility of $A_{1, m}$ 
by powers of $2$ and the
$3x+1$ problem. Namely $a_m := ord_2( A_{1, m}) - 1$ gives the number of iterations
of the $3x+1$ map $T(x)$ starting from $x_0=m$, in which the
parity of the iterates does not change, i.e
$$
m \equiv T(m) \equiv \cdots \equiv T^{a_{m}-1}(m) \not\equiv T^{a_m}(m) ~(\bmod~2).
$$
This is given as Theorem 6.1 of the paper.

\item
Paul  Andaloro (2000),
{\em On total stopping times under $3X+1$ iteration,}
Fibonacci Quarterly {\bf 38} (2000), 73--78. (MR 2000m:11024). \\
\newline 
\hspace*{.25in}
This paper shows various results on the minimal elements having a
given stopping time, where the ``stopping time'' is defined to
be the number of odd elements in the trajectory up to and
including $1$. It also obtains a new 
congruential ``sufficient set''  criterion to verify the
$3x+1$ Conjecture. It shows that knowing that the $3x+1$ Conjecture
is true for all $ n \equiv ~ 1~(\bmod~ 16)$ implies that it
is true in general.

\item
Paul Andaloro (2002),
{\em The $3X+1$ problem and directed graphs,}
Fibonacci Quarterly {\bf 40} (2002), 43--54. (MR 2003a:11018). \\
\newline 
\hspace*{.25in}
This paper considers various ``compressed'' versions of the
$3x+1$ graph, in which only a subset of the vertices
are retained
with certain directed  paths in original  $3x+1$ 
graph iterates of $T(\cdot)$ replaced by
single directed edges. The initial ``compressed'' graph corresponds to
odd integers, and the paper introduces
two further  ``compressed'' graphs with fewer allowed vertices.
In each case, the  $3x+1$ Conjecture is equivalent to 
the graph being weakly connected, i.e. being connected when
viewed as an undirected graph.
The paper shows that certain kinds of vertex pairs
in such graphs are weakly connected, typically for allowed 
vertices in certain
congruence classes $(\bmod~2^k)$ for small $k$.

\item
Stefan Andrei, Manfred Kudlek, Radu Stefan Niculescu  (2000),
{\em Some results on the Collatz problem},
Acta Informatica {\bf 37} (2000), 145--160.
(MR 2002c:11022). \\
\newline 
\hspace*{.25in}
This paper studies the internal structure of $3x+1$-trees.
Among other things, it  observes the that integers of the 
form $2^m k - 1$ eventually
iterate to the integers $3^m k - 1$, respectively.
and that integers  of the form
$2^{3m} k - 5$ iterate to the integers $3^{2m} k - 5$,
and integers  of the form $2^{11m} k - 17$ iterate to the
integers $3^{7m} k - 17$. These facts are related to the
cycles associated to $-1$, $-5$ and $-17$, respectively.


\item
Stefan Andrei, Wei-Ngan Chin, and Huu Hai Nguyen (2007+),
{\em A Functional View over the Collatz's Problem,}
preprint.\\
\newline 
\hspace*{.25in} 
This paper describes chains in the $3x+1$ tree using
a phrase structure grammar. 

\item
David Applegate and Jeffrey C. Lagarias (2003), 
{\em Lower bounds for the total stopping time  of $3x+1$ iterates,}
Math. Comp. {\bf 72} (2003), 1035--1049. (MR 2004a:11016).\\
\newline 
\hspace*{.25in} 
This paper proves  there are infinitely many positive $n$ which
have a finite total stopping time $\sigma_\infty(n) > 6.14316 \log n.$
It also shows that there is a positive $c$ such
that at least $ cx^{1/60}$ of all integers $1 < n \le x$
have a finite total stopping time $\sigma_\infty (n) > 5.9 \log n.$
The proofs are computer-intensive, and produce a ``certificate''
encoding a proof, which is based on 
a search of $3x+1$ trees
to depth $60$. The ``certificates'' are quite large, involving 
about $350$ million trees for the lower bound $6.14316\log n$,
which corresponds to a density of odd integers in
a trajectory (the'' ones-ratio'')  of $\frac {14}{29}\approx 0.483.$

This rigorous bound  is below the bound 
$\sigma_\infty(n) \approx 6.95212 \log n$ that one expects to
hold for almost all integers, which 
corresponds to a ones-ratio of $\frac{1}{2}$.
The paper gives heuristic arguments suggesting that the 
method of this paper might prove $\sigma_\infty(n) \approx 6.95212 \log n$ 
holds for infinitely many $n$,  but that it would likely require  a
search of  $3x+1$ trees to depth at least $76$. 
This would require a very large computation.

\item
David Applegate and Jeffrey C. Lagarias (2006), 
{\em The $3x+1$ semigroup}, 
J. Number Theory {\bf 177} (2006), 146--159. (MR 2006k:11037).\\
\newline 
\hspace*{.25in} 
This paper considers a weak version of the $3x+1$ problem
proposed by Farkas (2005). It considers the multiplicative
semigroup $\sR$ of positive rational numbers generated by
$\{ \frac{2n+1}{3n+2}: n \ge 0\}$ together with  $\{2\}$.
The {\em weak $3x+1$ conjecture} asserts that this semigroup 
contains all positive integers. The relation to the 
$3x+1$ problem is that the semigroup contains $1$ and
its generators encode the action of the inverse $3x+1$
function. It follows that the truth of the $3x+1$ conjecture
implies the truth of the weak $3x+1$ conjecture. 
This paper proves the
conjecture. Its main result shows that the semigroup $\sR$
consists of all positive rational numbers $\frac{a}{b}$ 
such that $3$ does not divide $b$. The proof is an induction
motivated by  certain results established in Lagarias (2006).

\item
Edward  Belaga and Maurice Mignotte (2000),
{\em Cyclic Structure of Dynamical Systems Associated
with $3x+d$ Extensions of Collatz Problem,}
U. Strasbourg report 2000-18, 57 pages. 
({\tt http://hal.archives-ouvertes.fr/IRMA-ACF}, file hal-00129656)\\
\newline 
\hspace*{.25in}
This paper is a theoretical and experimental study of
the distribution of number and lengths of finite cycles
to the $3x+ d$ map, for $d \equiv \pm 1~(\bmod 6).$
It contains much data and a wealth of interesting results
and conjectures.

\item
Edward Belaga (2003),
{\em Effective polynomial upper bounds to perigees and
numbers of $(3x+d)$-cycles of a given oddlength,}
Acta Arithmetica {\bf 106}, No. 2,  (2003), 197--206. (MR 2003m:11120). \\
\newline
\hspace*{.25in}
Let $d$ be a positive odd integer, and consider
the $3x+d$ map 
$T_d(x) = \frac{3x+d}{2}$ if $x$ is odd; $T_d(x) = \frac{x}{2}$
if $x$ is even, acting on the domain of positive integers.
This paper shows that for any cycle $C$ of the $3x+d$ map
of length $l$ containing $k$ odd elements, the
smallest element $prg(C)$  in the
cycle satisfies
$$
prg(C) \le \frac{d}{2^{l/k} - 3}.
$$  
From this follows 
$$ 
\log_2 3 < \frac{\mbox{length}(C)}{\mbox{oddlength}(C)} \le \log_2(d+3)
$$
The author shows that the upper bound is sharp, and gives
evidence that the 
lower bound is probably asymptotically approachable.
Using bounds from transcendence theory (linear forms in
logarithms) the author gives an upper bound  for  the
total number of cycles $U_{d,k}$ of odd-length $k$.
This  upper bound is  $d k^{c_0}$, for a constant $c_0$,
and states that  one may take $c_0=32$. He also shows that the
largest element $S$  of any such cycle is bounded above
by 
$$
S < d k^{c_0}(\frac{3}{2})^k.
$$

\item
Edward Belaga and Maurice Mignotte (2006a),
{\em Walking Cautiously into the Collatz 
Wilderness: Algorithmically,
Number Theoretically, and Randomly,}
Fourth Colloquium on Mathematics and Computer Science,
DMTS (Discrete Mathematics and Computer Science) Proceedings  {\bf AG}, 2006, 249--260.\\
\newline 
\hspace*{.25in}
This paper discusses many  open questions about the
$3x+1$ map and related maps, and recent new numerical
evidence supporting them. It formulates some new conjectures.

\item
Edward Belaga and Maurice Mignotte (2006b), 
{\em The Collatz problem and
its generalizations: Experimental Data. Table 1. Primitive cycles
of $3x+d$ mappings, } Univ. of Strasbourg preprint 2006-015, 9 pages+400+ page table .\\
\newline 
\hspace*{.25in}
This paper gives detailed tables of primitive cycles on the
positive integers for $3x+d$ maps 
for $1 < d < 20000$.  The authors conjecture they obtain the complete
list of such cycles, for these values of $d$. The text prior to the table formulates some 
interesting new
conjectures.

\item
Edward Belaga and Maurice Mignotte (2006c),
{\em The Collatz problem and
its generalizations: Experimental Data. Table 2. Factorization of
Collatz numbers $2^l-3^k$, }
Univ. of Strasbourg preprint 2006-018, 
6 pages text+156 page table. \\
\newline 
\hspace*{.25in}
These tables give known divisors of numbers of form $D=2^l-3^k$ for
$l \le 114$. Numbers $D$ of this form give $3x+D$ problems having
many primitive cycles.

\item
Vitaly Bergelson, Michal Misiurewicz and Samuel Senti (2006), 
{\em Affine Actions of a Free Semigroup on the Real Line,}
Ergodic Theory and Dynamical Systems {\bf 26} (2006), 1285--1305.
MR2266362 (2008f:37019). \\
\newline
\hspace*{.25in}
This extends the analysis of Misiurewicz and Rodrigues (2005)
to semigroups generalizing those associated to the $3x+1$ map. 
The current paper considers  the orbits of a
semigroup generated by $T_0(x)=ax$, $T_1(x)= bx+1$,
in which $0<a<1<b$, viewed as acting on the positive real
numbers $\RR_{\ge 0}.$ It defines various notions of
``uniform distribution'' on the positive real axis,
and derives various results concerning uniform distribution
of orbits of such semigroups. In particular, if the iterates of the transformation
are decreasing on average, and if one assumes the symbolic
dynamics of $T_0, T_1$ are drawn from an ergodic  invariant
measure $\nu$ on the shift space, then it proves the  existence and uniqueness of 
an  invariant 
measure $\mu$ on $\RR_{\ge 0},$
which depends on $\nu$. See Theorem F of the paper for a precise
statement. 

\item
Konstantin Borovkov and Dietmar Pfeifer (2000), 
{\em Estimates for the Syracuse problem via a probabilistic model,} 
Theory of Probability and its Applications {\bf 45}, No. 2 (2000), 300--310.
(MR 1 967 765). \\
\newline
\hspace*{.25in}
This paper studies a multiplicative random walk imitating the $3x+1$
iteration.
Let $X_0$ be given and set $Y_j = X_0 X_1 \cdots X_j$ where each
$X_i$ for $i \geq 1$ are i.i.d. random variables assuming the value
$\frac{1}{2}$ or $\frac{3}{2}$ with probability $\frac{1}{2}$ each.
Let $\sigma_\infty (X_0 , \omega )$ be a random variable equal to the
smallest $J$ such that $Y_j < 1$.
Then 
$E[ \sigma_\infty (X_0 , \omega ) ] = 
( \frac{1}{2} \log \frac{4}{3} )^{-1} \log X_0$ 
and the normalized variable
$$
\hat{\sigma}_\infty (X_0 , \omega ) = 
\frac{\sigma_\infty (X_0 , \omega ) - c_1 \log X_0}{c_2 ( \log X_0 )^{1/2}}
$$
with $c_1 = ( \frac{1}{2} \log \frac{4}{3} )^{-1} = 6.95212 $ and
$c_2 = c_1^{3/2} ( \frac{1}{2} \log 3)$, has a distribution converging to
a standard normal distribution as $X_0 \rightarrow \infty$ (Theorem 5).
Various refinements of this result are given.
Comparisons are made to empirical $3x+1$ data for
values around $n \approx 10^{6}$, which show good agreement.

\item
Barry Brent (2002+),
{\em $3X+1$ dynamics on rationals with fixed denominator,}
{\tt eprint: arXiv math.DS/0204170}. \\
\newline
\hspace*{.25in}
This paper reports on computer experiments looking for cycles 
with greatest common divior $1$ for the $3x+k$ problem, for various 
$k \equiv \pm 1 (\bmod~6)$. It suggests that for $k= 7$, $k=19$ and 
$k=31$ there is only one such cycle on the positive integers.
Various other statistics are reported on.

\item
Thomas Brox (2000), 
{\em Collatz cycles with few descents, }
Acta Arithmetica {\bf 92} (2000), 181--188. (MR 2001a:11032) \\
\newline
\hspace*{.25in}
This paper considers the variant of the $3x+1$ map, call
it $T_1$, that divides
out all powers of $2$ at each step. A cycle is written $\{x_1, x_2, ... x_m\}$
where each $x_i$ is an odd integer. A {\em descent} means
$|T_1(x)| < |x|$. Let $|C|$ denote the number of (odd) elements in
a $3x+1$ cycle $C$, and let $d(C)$ denote the number of descents
in $C$. 
The author proves that the the number of $3x+1$ cycles $C$ that
have $d(C) < 2 log |C|$ is finite. In particular the number of
Collatz cycles with $d(C) < r$ is finite for any fixed $r$. The
proof uses a bound of Baker and Feldman and is in principle effective.
Thus for each $r$ there exists an algorithm  to determine all cycles with
$d(C) < r$. This greatly strengthens the result of R. Steiner (1978)
who detemined all cycles with one descent.

\item
Paul S. Bruckman (2008),
{\em A proof of the Collatz conjecture,}
International Journal of Mathematical Education in Science and Technology,
{\bf 39}, No. 3 (2008), 403--407. [Erratum: {\bf 39}, No. 4 (2008), 567.] (MR 2009d:11043b) \\
\newline
\hspace*{.25in}
This paper asserts a proof of the Collatz conjecture. However the argument given
has a gap which leaves the proof incomplete.   The erratum points out this
gap and withdraws the proof. \\

The gap is, 
suppose $N_0$ is the starting value, and that $N_k$ is the $k$-th odd iterate to occur.
Let  $E_k$ denote the number of divisions by $2$ that  occur in reaching $N_k$, then
$ 2^{E_k} N_k - 3^k N_0 = S_k$, where $S_k$ is the positive
integer
$S_k = \sum_{j=0}^{k-1}  2^{E_j} 3^{k-1-j}.$  Next determine
 $(A_k, B_k)$ by requiring $2^{E_k} B_k - 3^k A_k = 1$, with $0 \le B_k < 3^k$.  
 The author notes that
there is an integer $T_k$ such that 
$$N_0= A_k+ T_k 2^{E_k},~~~~~ N_k = B_k + T_k 3^k.$$
Here $T_k$ depends
on $k$ and may be positive or negative. The author then argues by
contradiction, asserting  in Section 2 the claim
that the minimal counterexample $N_0$ must have $2^{E_k}< 3^k$ for 
all $k \ge 1$,  which would imply that the sequence of iterates of $N_0$ diverges. The
argument justifying this claim has a gap, it  supposes $2^{E_k} > 3^k$,
and asserts the contradiction that $N_k < N_0$. But $N_k \ge N_0$ 
may hold if $T_k$ is sufficiently negative, and in fact for every starting value
$n_0\ge 1$ the values $T_k \to - \infty$ as $k \to \infty$. This can be seen from
the equation (11) saying that
$n_0= A_k S_k + T_k 2^{N_k}$, noting that $n_0$ is fixed, $S_k \to \infty$
by its recurrence $S_{k+1} = 3S_k+ 2^{E_k}$ and $A_k \ge 1$.
(The author's argument as given would  prove the cycle starting at $1$ does not exist.)

\item
Charles C. Cadogan (2000),
{\em The $3x+1$ problem: towards a solution},
Caribbean J. Math. Comput. Sci. {\bf 10} (2000), paper 2, 11pp.
(MR 2005g:11032) \\
\newline
\hspace*{.25in}
The paper studies trajectories  of the $3x+1$ problem.
calling two integers $n_1$ and $n_2$ equivalent,
written $n_1 \sim n_2$, if their trajectories eventually
coalesce. Various results are obtained  giving sufficient
conditions for equivalence. The author conjectures that
for each positive odd integer $n$, $9n+4 \sim 3n+1$.  His main result is
that this conjecture implies the truth of the $3x+1$ Conjecture.

\item
Charles C. Cadogan (2003),
{\em Trajectories in the $3x+1$ problem,}
J. of Combinatorial Mathematics and Combinatorial Computing,
{\bf 44} (2003), 177--187. (MR 2004a:11017)  \\
\newline
\hspace*{.25in}
This paper describes  various pairs of  trajectories
that coalesce under the $3x+1$ iteration. For example
the trajectories of $3n+1$ and $4n+1$ coalesce,
and the trajectory of $16k + 13$ coalesces with 
that of $3k+4$. The main result (Theorem 3.9) gives a certain
infinite family of coalescences. 

\item
Charles C. Cadogan (2006),
{\em A Solution to the $3x+1$ Problem,}
Caribbean J. Math. Comp. Sci. {\bf 13} (2006), 1--11. \\
\newline
\hspace*{.25in}
This paper asserts a proof of the $3x+1$ conjecture. However 
the argument given
has a gap which leaves the proof incomplete.
Namely,  on the line just before equation (2.6) the
expression $1+ 2 t_{i, j} \sim 1+ 3n_{i+1, j}$ should instead read
$1+ 2t_{i,j} \sim 1+ 3n_{i, j}$, as given by equation (2.5). 
Hence instead of obtaining equation (2.6) in
the form  $t_{i,j} \sim t_{i+1, j} \sim t_{i+2, j}$,
one only obtains $t_{i, j} \sim t_{i+1,j}$.
This renders the  proof of Theorem 2.15  incomplete, as it 
depends on equation (2.6).  Then the induction step in Lemma 3.1
cannot be completed, as it depends on Theorem 2.15. Finally 
the main result Theorem 3.3 has a gap since it depends on
 Lemma 3.1.

\item
M\'{o}nica del Pilar Canales Chac\'{o}n 
and  Michael Vielhaber (2004),
{\em Structural and Computational Complexity of Isometries
and their Shift Commutators,} Electronic Colloquium
on Computational Complexity, Report No. 57 (2004), 24 pp. (electronic). \\
\newline
\hspace*{.25in}
The paper considers functions on $f: \{0, 1\}^{\infty} \to 
 \{0, 1\}^{\infty}$ 
computable by invertible transducers. They give several
formulations for computing such maps and consider 
several measures of 
computational complexity of such functions, including
tree complexity $T(f, h)$, which measures the
local branching of a tree computation, where $h$ is the tree height
of a vertex. They also study the bit complexity $B(f,n)$ which is the
complexity of computing the first $n$ input/ output symbols.
Tree complexity is introduced in
H. Niederreiter and M. Vielhaber, J. Complexity {\bf 12} (1996), 
187--198 (MR 97g:94025).

The $3x+1$ function is considered as an example
showing that some of the general  complexity bounds 
they obtain are sharp.
Interpreting the domain $\{0, 1\}^{\infty}$
as the $2$-adic integers, the map $Q_{\infty}$
associated to the $3x+1$ map
given in Lagarias (1985) [Theorem L] is a function of this kind.
It is invertible and the inverse map $Q_{\infty}^{-1}$ is
studied in Bernstein (1994)
and Bernstein  and Lagarias (1996). In Theorem 33
the authors give a 5-state shift automaton that computes
the ``shift commutator'' of the $3x+1$ function, which
they show takes a $2$-adic integer $a$ to $a$ if $a$ is even, and to
$3a+2$ if $a$ is odd. 
In Theorem 34 they deduce that the tree complexity of
$Q_{\infty}$ is bounded by a constant. 
Here $Q_{\infty}$ corresponds to their function $T({\bf c}, \cdot)$.


\item
Mark Chamberland (2003),
{\em Una actualizachio del problema $3x+1$,}
Butletti de la Societat Catalana, {\bf 22} (2003) 19--45.
 (MR 2004i:11019).\\
\newline
\hspace*{.25in}
This is a survey paper (in Catalan) describing  recent results on
the $3x+1$ problem, classified by area. 

{\em Note.} An English version of this paper: ``An Update on the 
$3x+1$ Problem'' is posted on
the author's webpage: {\tt http://www.math.grin.edu/~chamberl/}

\item
Dean Clark, 
Periodic solutions of arbitrary length in a simple integer iteration,
Advances in Difference Equations, {\bf 2006} Article 35847, pp. 1--9.\\
\newline
\hspace*{.25in}
This paper studies the second order nonlinear recurrence
$$
y_{n+1} = \lceil a y_n\rceil - y_{n-1}
$$
in which the initial conditions $(y_0, y_1)$ are integers, and $a$
is a constant with 
$\{ a\in \RR: |a|<2, \}$. These generate an infinite
sequence. For  $a= \frac{3}{2}$,  presuming integer initial conditions,
the recurrence can be rewritten as
$$
y_{n+1} :=\left\{
\begin{array}{cl}
\df{3y_n +1}{2} - y_{n-1}  & \mbox{if} ~~y_n\equiv 1~~ (\bmod ~ 2) \\
~~~ \\
\df{3y_{y}}{2} - y_{n-1}  & \mbox{if}~~ y_n \equiv 0~~ (\bmod ~2 )~.
\end{array}
\right.
$$
which the author views as a second-order
analogue of the $3x+1$ iteration. The main result of the paper 
shows that for fixed $|a|<2$, all orbits of the recurrence are purely
periodic. 

In addition, the author shows when $a= \frac{p}{q}$
is a non-integer rational with $|a|<2$, then  there are orbits of arbitrarily large period.

The author presents some  computer plots of values $T(y_{n-1}, y_n) =(y_n, y_{n+1})$ in the plane, 
For $a = \frac{1+ \sqrt{5}}{2}$ and other values, the periodic orbits plotted this way can
appear complex, approximating fractal-like shapes.

\item
Lisbeth De Mol (2008),
{\em Tag Systems and Collatz-like functions},
Theoretical Computer Science {\bf 390} (2008), 92--101.\\
\newline
\hspace*{.25in}
A Post tag stystem $\sT$ with alphabet size $\mu$ and shift $\nu$ consists  of 
an alphabet $\sA= \{a_1, a_2, \cdots , a_{\mu}\}$ and rules $a_i \mapsto E_i$,
where $E_i$ is a finite word with symbols drawn  from this alphabet. Given an input word $w$
with leftmost symbol $a_i$, the tag system iteration produces
an output word $w'$ which chops off the $\nu$ leftmost
symbols and tags on to  the right  end  of $w$ the word $E_i$,  thus 
the new word $w'$ has length $|w'| = |w| +|E_i|-\nu.$
 Given an initial word $w_0$, the tag system produces subsequent words
 $w_1, w_2, ...$. It is said to halt at step $n$ if $w_n$ is  the empty word.
The {\em halting problem} for a tag system $\sT$ is to determine for each possible input word $w_0$
whether or not the tag system eventually halts on that input. The {\em reachability problem} for $\sT$ is,
given a word $ x$,  to determine for each input $w_0$ whether or not the word $x$
is  eventually reached under repeated iteration. 
Let $TS(\mu, \nu)$ denote the collection of all Tag systems with parameters $(\mu, \nu)$.
In 1921 Emil Post showed the halting and reachability  problems
are decidable for all tag systems in $TS(2,2)$; his proofs were not published.  
He could not resolve the case of $TS(2,3).$
In 1961 Minsky  showed that the 
halting problem, and hence the reachability problem are undecidable for tag systems
in general. It was  later shown that a universal Turing machine is encodable as a tag
system in  $TS(576, 2)$, so both the halting and reachability  problems are unsolvable in this class. 

In this paper the author shows the $3x+1$ problem is encodable as a reachability
problem in $TS(3,2)$. The tag system takes $\sA=\{ a, b,c\}$ with rules
$a \mapsto bc, ~b \mapsto a, ~c \mapsto aaa$. The reachability question concerns  reaching $x= a$, Starting
from $w_0= a^n$ (word repeating the letter $a$ $n$ times)  the next word that
is a power of $a$ that the tag system reaches is
 $a^{T(n)}$, where $T(n)$ is the $3x+1$ function. To answer this reachability problem
 for all iterates requires solving the $3x+1$ problem. It follows that in should be a difficult 
 problem to find
 a decision procedure for the reachability problem on the class $TS(3,2)$.
 The author also shows that the iteration of a generalized Collatz function
 which is affine on congruence classes $(\bmod~d)$ is encodable in a tag system
 in some $TS(\mu, d)$ with $\mu \le 2d+3$.
 
 \item
 Diego Domenici (2009)
 {\em A few observations on the Collatz problem,}
 Inter. J. Appl. Math. Stat. {\bf 14} (2009), No. J09, 97--107.
 (MR 2524878)
 \newline
\hspace*{.25in}
The paper proves results encoding the $3x+1$ conjecture as a problem of representing
integers in particular forms. It shows a necessary and sufficient condition for the
$3x+1$ conjecture to hold is that each positive integer $n$ has a representation
$$
n = \frac{2^m}{3^l} - \sum_{k=1}^1 \frac{2^{b_k}}{3^k}
$$
with $1 \le l \le m-3$ and
$0 \le b_1 < b_2 < \cdots < b_l \le  m-4.$
 
 \item
Jean-Guillaume Dumas (2008),
{\em Caract\'{e}risation des quenines et leur repr\'{e}sentation spirale}
Math\'{e}matiques et Sciences Humaines
[Mathematics and Social Science], {\bf 184} No. 4, 9--23.\\
\newline
\hspace*{.25in}
This paper solves the problem raised by Queneau (1963) on
allowable spiral rhyme patterns generalizing the `sestina" rhyme pattern 
of Arnaut Daniel, a 12-th century troubadour. 
One considers iterations of the $(3x+1)$-like function 
$$
\sigma_n(x) := \left\{ 
 \begin{array}{cl}
\frac{x}{2} & ~~\mbox{if} ~~x  ~~\mbox{is ~even}\\
~&~\\
\frac{2n+ 1-x}{2} & ~~\mbox{if} ~~x  ~~\mbox{is ~odd}\\
\end{array}
\right.
$$
on the domain $\{1, 2, ..., n\}$. This  function gives a permutation in
the symmetric group $S_n$ 
of this domain, which is called by the author a {\em quenine}.  We also
recall the inverse permutation
$\delta_{2,n}(x) = \sigma_n^{-1}(x)$ given by 
$$
\delta_{2,n}(x)  := \left\{ 
 \begin{array}{cl}
2x & ~~\mbox{if} ~~1 \le x \le \frac{n}{2}\\
~&~\\
2n+1-2x & ~~\mbox{if} ~~\frac{n}{2} < x \le n.\\
\end{array}
\right.
$$
The Raymond Queneau numbers $n$ (or $2$-admissible numbers) are those numbers for which 
this permutation is a cyclic permutation.

The author first recalls (Theorem 1) the results of Berger (1969), which
used the inverse permutation $\delta_{2,n}(x)$ above to
put restrictions on Queneau numbers. 
The author then proves (Theorem 2) that an integer $n$ is a Queneau
number  if and only if either $p=2n+1$ is prime, and either
(i) $p \equiv 3$ or $5~(\bmod~8)$ and $2$ is a primitive root of $p$,
i.e. $ord_p(2)= 2n$, or (ii) $p \equiv 7~(\bmod~8)$ and $ord_p(2) = n$.
This establishes a corrected form of a Conjecture of Roubaud (1993).

The author goes on to consider generalizations of quenines
introduced by Roubaud (1993). He gives pictures showing the various
spiral patterns associated with such permutations. To start, he sets 
$$
\delta_{3,n}(x) := \left\{ 
 \begin{array}{cl}
3x & ~~\mbox{if} ~~1 \le x \le \frac{n}{3}\\
~&~\\
2n+1-3x & ~~\mbox{if} ~~\frac{n}{3} < x \le \frac{2n}{3}\\
~&~\\
3x- (2n+1) & ~~\mbox{if} ~~\frac{2n}{3} < x \le n.\\
\end{array}
\right.
$$
More generally, one can replace $3$ with any $g \le n$, to
get a $g$-spiral permutation, whose definition can be worked out
from a picture of the spiral. 
Call such a permutation {\em $g$-admissible}
if it is a cyclic permutation. Theorem 3 gives a necessary
and sufficient condition for $n$ to be $g$-admissible, which is
that $p=2n+1$ is prime and either (i) $g$ is a primitive root $(\bmod ~p)$,
or (ii) $n$ is odd, and $ord_p(g) = n$.
The paper includes a table giving for $n \le 1000$ such that $p=2n+1$ is
prime the minimal $g$ for which $n$ is $g$-admissible. 

Various other types of permutations are considered: a {\em p\'{e}recquine}
is one associated to the permutation
$$
\pi_n(x) := \left\{ 
 \begin{array}{cl}
2x & ~~\mbox{if} ~~1 \le x \le \frac{n}{2}\\
~&~\\
2x-(n+1)& ~~\mbox{if} ~~\frac{n}{2} < x \le n.\\
\end{array}
\right.
$$
This pattern  is named after Georges Perec, an Oulipo member, cf. Queneau (1963).  
Theorem 5 shows that this permutation is cycle if and only if $p=n+1$ is prime and 
$2$ is a primitive root $(\bmod ~p)$.

\item
Jeffrey P. Dumont and Clifford A. Reiter (2001),
{\em Visualizing Generalized 3x+1 Function Dynamics,}
Computers and Graphics {\bf 25} (2001), 883--898. \\
\newline
\hspace*{.25in}
This paper describes numerical and graphical experiments
iterating generalizations of the $3x+1 $ function.
It plots  basins of attraction and
false color pictures of escape times for various generalizations of the
$3x+1$ function to the real line, as in Chamberland (1996),
and to the complex plane, as in Letherman, Schleicher and Wood (1999).
It introduces a new generalization to the complex
plane , the {\em winding 3x+1 function},
$$W(z) := \frac{1}{2}\left( 3^{ mod_2 (z)} z  + mod_2(z) \right),$$
in which
$$ mod_2(z) := \frac{1}{2}( 1 - e^{\pi i z}) 
= (\sin \frac{\pi z}{2})^2 - \frac{i}{2} \sin \pi z,$$
Plots of complex basins of attraction for Chamberland's function  appear
to have a structure resembling the Mandlebrot set, while the 
basins of attraction of the winding
$3x+1$ function seems to have a rather different structure.
The programs were written in the computer language J.

\item
Jeffrey P. Dumont and Clifford A. Reiter (2003),
{\em Real dynamics of a $3$-power extension of the $3x+1$ function,}
Dynamics of Continuous, Discrete and Impulsive Systems, Series A:
Mathematical Analysis 
{\bf 10} (2003), 875--893. (MR2005e:37099).\\
\newline
\hspace*{.25in}
This paper studies the real dynamics of the function
$$T(x) := \frac{1}{2}(3^{mod_{2}(x)}x + mod_{2}(x)),$$
 in
which the function
$$
mod_{2}(x) := (\sin\frac{1}{2}\pi x)^2.
$$
This function agrees with the $3x+1$-function on the
integers. The authors show that this function
has negative Schwartzian derivative on the
region $x > 0$. They study its periodic
orbits and critical points, and
show that any cycle of positive integers
is attractive. They define an extension of the notion of
total stopping time to all real
numbers $x$ that are attracted to the
periodic orbit $\{1,2\}$, representing the
number of steps till the orbit enters the
immediate basin of attraction of this attracting
periodic orbit. 
They formulate the
{\em odd critical point conjecture}, which asserts for
an odd positive integer $n \ge 3$ with associated nearby
critical point $c_n$, that the critical point $c_n$ is attracted
to the periodic orbit $\{1, 2\}$, and that $c_n$ and $n$ have
the same total stopping time.

\item
Hershel M. Farkas (2005),
{\em Variants of the $3N+1$ problem and multiplicative semigroups}, 
In:  Geometry, Spectral Theory, Groups and Dynamics: Proceedings in
Memory of Robert Brooks, 
Contemporary Math., Volume 387, Amer. Math. Soc., Providence, 2005, pp. 121--127.
(MR 2006g:11052)\\
\newline
\hspace*{.25in}
This paper formulates some weakenings of the $3x+1$ problem
where stronger results can be proved. It first shows
that iteration of the map
$$
F(n) = \left\{
\begin{array}{cl}
 \frac{n}{3} &  \mbox{if}~~n \equiv~~~~0 ~(\bmod~ 3));
~~~ \\
\frac{3n+1}{2} & \mbox{if}~~ n \equiv ~7 ~\mbox{or}~ 11~(\bmod~12);
~~~\\
\frac{n+1}{2} & \mbox{if}~~ n\equiv~1~ \mbox{or}~ 5~(\bmod~12);
\end{array}
\right.
$$
on the positive integers has all trajectories 
get to $1$.
The trajectories of the iterates of these functions
can have arbitrarily long subsequences on which the
iterates increase. 
The author then asks questions of the type:
``Which integers can be
represented in a multiplicative semigroup whose
generators are a specified infinite set of rational numbers?''
He proves that the integers represented
by the multiplicative semigroup
generated by  $\{ \frac{d(n)}{n}:~ n \ge 1\}$,
where $d(n)$ is the divisor function, represents
exactly the set of positive odd integers. The analysis
involves the function $F(n)$ above. Finally the author
proposes as an open problem a  weakened version of
the $3x+1$ problem, which asks:
``Which integers are represented by the
multiplicative semigroup generated by
$\{ \frac{2n+1}{3n+2}: n \ge 1\}$ together with $\{2\}$?''
The truth of the $3x+1$ Conjecture implies that
all positive integers are so represented. 

{\em Note.} Applegate and Lagarias (2006)
prove that all positive integers are represented
in the semigroup above. 


Of the author's three results, 
Theorem 1 holds for $\tilde{T}(x)$. The arguments of Theorems 2 and 3
presented appear to apply to $\tilde{T}(x)$  as well; line 4 of Theorem 3 needs to be
modified to $\tilde{T}^{(L+1)}(n) \equiv \tilde{T}^{(L)}(n)+ 1~(\bmod~2)$. 
The ``proof" of Theorem 2 seems  not  to adequately relate the mathematics
and the metamathematics. 

\item
David Gluck and Brian D. Taylor (2002),
{\em A new statistic for the $3x+1$ problem,}
Proc. Amer. Math. Soc. {\bf 130} (2002), 1293--1301.
(MR 2002k:11031). \\
\newline
\hspace*{.25in}
This paper considers iterations of the Collatz function 
$C(x)$. If $\ba = (a_1, a_2, ..., a_n)$ is a finite Collatz
trajectory starting from $a_1$, with $a_n = 1$ being the
first time $1$ is reached, they assign the statistic
$$C(\ba) = \frac{a_1a_2 + a_2 a_3 +... + a_{n-1}a_n + a_n a_1}
{a_1^2 + a_2^2+...+ a_n^2}.$$ 
They prove that $\frac{9}{13} < C(\ba) < \frac{5}{7}$.
They find sequences of starting values that approach the upper
and lower bounds, given that the starting values terminate.

\item
Jeffrey R.  Goodwin (2003),
{\em Results on the Collatz conjecture,} 
Annalele Stiintifice ale Universitatii ``Al. I. Cuza'' din lasi serie noua.
Informatica  (Romanian), {\bf XIII} (2003) pp. 1--16. 
MR2067520 (2005b:11025).
[Scientific Annals of the ``Al. I. Cuza'' University of Iasi,
Computer Science Section, Tome XIII, 2003, 1--16]\\
\newline
\hspace*{.25in}
This paper partitions the inverse iterates of $1$ under
the $3x+1$ map into various subsets, and studies their
internal recursive structure.

\item
He, Sheng Wang (2003),
{\em $3n+1$ problem's simplifying and structure property of $\{T(n)\}$$(n \in \NN)$}
(Chinese),
Acta Scieniarum Naturalium Universitatis Neimonggol 
[Nei Menggu da xue xue bao. Zi ran ke xue] (2003), No. 2. \\
\newline
\hspace*{.25in}
[I have not seen this paper.]

\item
Kenneth  Hicks, Gary L. Mullen, Joseph L. Yucas and Ryan Zavislak (2008),
{\em A Polynomial Analogue of the $3N+1$ Problem?}, 
American Math. Monthly {\bf   115}  (2008), No. 7, 615--622.\\
\newline
\hspace*{.25in}
This paper studies an iteration on polynomials resembling
the $3x+1$ problem. It considers the function on 
the polynomial ring $GF_2[x]$, where $GF_2$ is the finite field
with $2$ elements,  given by 
\begin{equation} ~\label{eq1}
C_{1}(f(x))  = \left\{
\begin{array}{cc} (x+1) f(x) + 1 & \mbox{if}~ f(0) \ne 0 \\
     \frac{f(x)}{x} & \mbox{if} ~ f(0) = 0. \end{array} 
     \right.
 \end{equation}
They show  that the iteration always converges
to a constant; for  an $n$-th degree starting polynomial in
at most $n^2+2n$ steps. They also show there exist starting polynomials
of degree $n$ which take at least $3n$ steps to get to a constant. 
They observe that their result applies more generally
to  the map  $C_F$ acting on monic polynomials
$f(x) = x^n +a_{n-1} x^{n-1} + \cdots + a_1 x + a_0 \in F[x]$, where $F$ is an arbitrary field,   by:
\begin{equation} ~\label{eq1}
C_{F}(f(x))  = \left\{
\begin{array}{cc} (x-\frac{a_0}{a_j}) f(x) + \frac{a_0^2}{a_j} & \mbox{if}~ f(0) \ne 0 \\
     \frac{f(x)}{x} & \mbox{if} ~ f(0) = 0, \end{array} 
     \right.
 \end{equation}
 where 
 $a_j$ is the smallest value $1 \le j \le n$ such that $a_j \ne 0$.
 Here a degree $n$ polynomial takes at most $n^2+2n$ steps to
 get to a constant polynomial.

{\em Note.} Iterations of a polynomial ring over a finite field were first
introduced in Matthews and Leigh (1987). They exhibited a mapping
having a provably divergent trajectory.

\item
Wernt Hotzel (2003),
{\em Beit\"{a}ge zum $3n+1$-Problem,}
Dissertation: Univsit\"{a}t Hamburg 2003, 62pp. [Zbl 1066.11501]\\
\newline
\hspace*{.25in}
This thesis consists of 8 short chapters. 

Chapter 2 describes symbolic dynamics of the $3x+1$ map $T$,
calling it  the Collatzfunction $d$, and noting it extends to the
domain of 2-adic integers $\ZZ_2$. He encodes a symbolic
dynamics of forward iterates, using binary
labels (called $I$ and $O$, described here as $1$ and $0$).
He introduces the $3x+1$ conjugacy map, denoting it $T$.

Chapter 3 studies the set of rational periodic points,  letting $\QQ_2$
denote the set of rationals with odd denominators. He observes that there
are none in the open interval $(-1, 0)$. 

Chapter 4 studies periodic cycles using Farey sequences, following
Halbeisen and Hungerb\"{u}hler (1997).

Chapter 5 studies unbounded orbits, observing that for rational inputs
$r \in \QQ_2$ has an unbounded orbit if and only if $r \not\in \QQ_2$.

Chapter 6 studies rational  cycles with a fixed denominator $N$,

Chapter 7 studies periodic points of the $3x+1$ conjugacy map.
The fixed point $\frac{1}{3}$ is known. An algorithm which searches
for periodic points is described.

Chapter 8 describes an encryption algorithms based on the $3x+1$ function.

\item
Huang, Guo Lin  and Wu, Jia Bang (2000),
{\em One-to-one correspondence between the natural numbers and the
parity vectors in the Collatz problem} (Chinese),
J. of South-Central University for the Nationalities, Natural Science Ed.
[Zhong nan min zu  da xue xue bao. Zi ran ke xue   ban] {\bf 19} (2000), No. 3, 59--61.\\
\newline
\hspace*{.25in}
English summary: "It is shown in this paper that $M_N^{\circ}= \{ 0,1,2, \cdots, 2^N-1\}$
and $V_N$ denotes the set of all $v_n$ of truncations up to the $N$-th term,
 viz. $\{ x_0, x_1, ..., x_{N-1}\}$
of the parity vector $v= \{ x_0, x_1, \cdots\}$, if $m \in M_N^{\circ}$ and $m \rightarrow v_N(m)$
then the mapping $\sigma: M_N^{\circ} \to V_N$ is one-to-one. The following lemmas are
from this. Lemma 1. Let $M_N=\{ 1, 2, \cdots, 2^N\}$, if $m \in M_N$ and $m \rightarrow v_N(m)$
then the correspondence is one-to-one. Lemma 2. Let ${\bf v}$ denote the set of all
parity vectors ${\bf v}=\{ x_0, x_1, \cdots\}$, if $m \in \NN$ and $m \to {\bf v}(m)$ then the correspondence
is one-to-one. Thus the investigation on any natural number can be converted to the investigation
of its parity vector."

{\em Note.} This result is implicit in  Terras (1976), and properties of
the parity vector are described in Lagarias (1985).

\item
Yasuaki Ito and Koji Nakano (2009),
{\em A Hardware-Software Cooperative Approach to the Exhaustive verification of the
Collatz conjecture,}
International Symp. on Parallel and Distributed Processing with Applications,
2009, pp. 63--70.\\
\newline
\hspace*{.25in}
[I have not seen this paper.]

\item
Ke, Wei (2001),
{\em Generalization Concerning the $3x+1$ Problem} (Chinese),
Journal of Jinan University (Science and Technology) 
[Jinan da xue xue bao. Zi ran ke xue ban]  (2001), No. 4.\\
\newline
\hspace*{.25in}
[I have not seen this paper.]

\item
Ke, Yong-Sheng (2005),
{\em The proof of the hypothesis about $"3x+1"$} (Chinese),
Journal of Tianjin Vocational Institute [Teacher's College]
[Tianjin zhi ye ji shu shi fax xue yuan xue bao] (2005), No. 2.\\
\newline
\hspace*{.25in}
[I have not seen this paper.]

\item
Immo O. Kerner (2000),
{\em Die Collatz-Ulam-Kombination CUK},
Rostocker Informatik-Berichte No. 24 (2000), preprint. \\
\newline
\hspace*{.25in}
This paper attibutes the Collatz problem to Collatz in 1937.
Let $c(n)$ be a function counting the number of iterations of
the Collatz function $C(n)$ to get to $1$, with $c(1) = 0$
and $c(3) = 7$, and so on. It reviews some basic results on the
Collatz function, and introduces an auxiliary two-variable
function $B(n, t)$ to describe the iteration, and plots some
graphs of its behavior.

\item
Stefan Kohl (2005),
{\em Restklassenweise affine Gruppen}, Universit\"{a}t Stuttgart,
Ph. D. Dissertation, 2005.
(eprint: {\tt http://deposit.ddb.de/dokserv?idn=977164071})
\newline
\hspace*{.25in}

In this thesis, the author studies the semigroup
$Rcwa(\ZZ)$consisting of all functions $f: \ZZ \to \ZZ$
for which there exists a modulus
 $m=m(f)$ such that the restriction of $f$
to each residue class $(\bmod~m)$ is an  affine map.
It also considers the group $RCWA(\ZZ)$ 
consisting of the set of invertible
elements of $Rcwa(\ZZ)$. The group $RCWA(\ZZ)$ is
a subgroup of the 
 infinite permutation
group of the  the integers. Both the $3x+1$ function
and the Collatz function belong to the semigroup $Rcwa(\ZZ)$. 
The original Collatz map (see Klamkin (1963))
given by $f(3n)=2n$, $f(3n-1)=4n-1$, $f(3n-2)=4n-3$
is a permutation belonging to $RCWA(\ZZ)$. 

Some of the results of the thesis are
as follows.  The group $RCWA(\ZZ)$ is
not finitely generated (Theorem 2.1.1). It
has finite subgroups of any isomorphism type (Theorem 2.1.2).
It has  a trivial center (Theorem 2.1.3). It 
acts highly transitively on $\ZZ$ (Theorem 2.1.5).
All nontrivial normal subgroups also act highly transitively
on $\ZZ$, so that it has no notrivial solvable
normal subgroup (Corollary 2.1.6). 
It has an epimorphism $sgn$ onto $\{ 1, -1\}$, so has a
normal subgroup of index $2$ (Theorem 2.12.8).
Given any two subgroups, it has another subgroup
isomorphic to their direct product(Corollary 2.3.3).
It has only finitely many conjugacy classes of  elements
having a given odd order, but it has  infinitely many conjugacy classes
having any given even order (Conclusion 2.7.2). 

The author notes that the $3x+1$ function can be 
embedded as a permutation in $RCWA(\ZZ \times \ZZ)$,
as $(x,y) \mapsto (\frac{3x+1}{2}, 2y)$ if $x \equiv 1(\bmod~2)$;
$\mapsto (\frac{x}{2}, y)$ if $x\equiv 0, 2 (\bmod~6)$;
$\mapsto  (\frac{x}{2}, 2y+1)$if $x\equiv 4 (\bmod~6)$,
where it represents the iteration projected onto the $x$-coordinate.

The  author develops an algorithm  for efficiently
computing periodically linear functions, whether they
are permutations or not.
Periodically linear functions are functions which are
defined as affine functions on each congruence class
$j ~(\bmod~M)$ for a fixed modulus $M$, as 
in Lagarias (1985). The $3x+1$ function 
is an example of such a function.
The author has written a corresponding package RCWA (Residue Class-Wise
Affine Groups) for  the  computational algebra and group theory 
system GAP (groups, Algorithms, Programming).
This package is available for download at: \\
~~~~~~~~{\tt http://www.gap-system.org/Packages/rcwa.html}. \\
A manual for RCWA 
can  be   found at: \\
~~~~~~~~{\tt http://www.gap-system.org/Manuals/pkg/rcwa/doc/manual.pdf}

\item
Stefan Kohl (2007),
{\em Wildness of iteration of certain residue class-wise affine mappings},
Advances in Applied Math. {\bf 39} (2007), 322--328. (MR 2008g:11041)\\
\newline
\hspace*{.25in}
A mapping $f: \ZZ \to \ZZ$ is called {\em residue class-wise affine} (abbreviated
RCWA)  if it is affine
on residue classes $(\bmod ~m)$ for some fixed $m \ge 1$.
(This class of functions was termed {\em periodically linear}
in Lagarias (1985).) The smallest such $m$ is called the {\em modulus}
of $f$.  This class of functions is closed under pointwise addition
 and under composition.
A function $f$ is called {\em tame} if the modulus of its $k$-th iterate
remains bounded as $k \to \infty$; it is {\em wild} otherwise. The
author shows that if $f: \ZZ \to \ZZ$ is an RCWA -function,
which is surjective, but not injective,
then $f$ is necessarily wild. The paper also presents counterexamples showing
that each of  the three other possible
combinations of  hypotheses of (non-)surjectivitiy or of (non)-injectivity 
of $f$ permits no conclusion whether it is tame or wild.

\item
Stefan Kohl (2008a),
{\em On conjugates of Collatz-type mappings,}
Int.  J. Number Theory {\bf 4}, No. 1 (2008), 117--120. (MR 2008m: 11055) \\
\newline
\hspace*{.25in}
A map $f: \ZZ \to \ZZ$ is said to be {\em almost contracting} 
if there is a finite set $S$ such that
every trajectory of $f$ visits this finite set. This property holds if
and only if there is a permutation $\sigma$ of the integers such
that $g=\sigma^{-1} \circ  f \circ \sigma$ decreases absolute value off a finite set,
a property that is called {\\em monotonizable}.
Suppose that $f: \ZZ \to \ZZ$ is a surjective, but not injective, RCWA mapping
(see Kohl(2006+) for a definition)
 having the property that  the preimage set of any integer under
$f$ is finite. The main result asserts that if $f$ is almost contracting and
$k$ such that the $k$-th iterate $f^{(k)}$ decreases almost all integers,
then any permutation $\sigma$ that establishes the almost contracting
property of $f$ cannot itself be an RCWA mapping. 

 The $3x+1$ function
$T(x)$ is believed to be a function satisfying the hypotheses of the author's main
result: It is surjective but not injective, and is believed to 
be almost contracting.  The almost contacting property for
the $3x+1$ map is equivalent to  establishing 
that its iteration on $\ZZ$ has only finitely many cycles  and no divergent trajectories.

\item
Stefan Kohl (2008b),
{\em Algorithms for a class of infinite permutation groups,}
J. of Symbolic Computation {\bf 43}, No. 8 (2008), 545--581. (MR 2009e:20003) \\
\newline
\hspace*{.25in}
A mapping $f: \ZZ \to \ZZ$ is called {\em residue-class-wise  afffine}
(RCWA) if there is a positive integer $m$ such that it is an affine
mapping when restricted to each residue class $(\bmod~m)$.
The $3x+1$ mapping $T$ is of this kind, as is a permutation
constructed by Collatz.
This paper describes a collection of algorithms and methods
for computing in permutation groups and monoids whose members are all
 RCWA mappings.

\item
Stefan Kohl (2010),
{\em A simple group generated by involutions interchanging residue classes
of the integers}, Math. Z. {\bf 264} (2010), no. 4, 927--938.\\
\newline
\hspace*{.25in}
The author defines a simple group of permutations of the integers, generated
by compositions of certain RCWA functions.
\item
Pavlos B.  Konstadinidis (2006),
{\em The real $3x+1$ problem},
Acta Arithmetica {\bf 122} (2006), 35--44. (MR 2007c:11029)\\
\newline
\hspace*{.25in}
The author extends the $3x+1$ function to the
real line as:
$$
U(x) = \left\{
\begin{array}{cl}
\df{3x+1}{2} & \mbox{if} ~~\lfloor x \rfloor \equiv 1~~ (\bmod~2) ~. \\
~~~ \\
\df{x}{2} & \mbox{if}  ~~\lfloor x \rfloor  \equiv 0~~ (\bmod~2)~.
\end{array}
\right.
$$
The paper shows
that the only periodic orbits of the function $U(x)$ on
the positive real numbers are those on the positive
integers. The paper also considers some related functions.

\item
Alex V. Kontorovich and Steven J. Miller (2005),
{\em Benford's law, values of $L$-functions, and the
$3x+1$ problem}, Acta Arithmetica {\bf 120} (2005), 269--297. 
(MR 2007c:11085). \\
\newline
\hspace*{.25in}
Benford's law says that the leading digit of decimal
expansions of certain sequences are not uniformly
distributed, but have the probability of digit $j$
being $\log_{10}(1+\frac{1}{j})$. This paper gives
a general method for verifying Benford's law for
certain sequences. These include special values of
L-functions and ratios of certain $3x+1$ iterates,
the latter case being covered in Theorem 5.3.
It considers the $3x+1$ iteration in the form of Sinai(2003a)
and Kontorovich and Sinai (2002), 
which for an odd integer $x$ has $m(x)$ being
the next odd integer occurring 
in the $3x+1$ iteration. For a given real base $B> 1$ it looks at
the distribution of the 
quantities $\log_B ( x_m/ (\frac{3}{4})^m x_0) (\bmod~1)$
as $x_0$ varies over odd integers in  $[1, X]$ and $X \to \infty$.
It then takes a second  limit as $m \to \infty$
and concludes that the uniform distribution is approached,
provided  $B$ is such that $\log_2 B$ is
an irrational number of finite Diophantine type.
The case $B=10$ corresponds to Benford's law. The
theorem applies when $B=10$ because  $\log_2 10$
is known to be of finite
Diophantine type by A. Baker's results on linear forms
in logarithms. A main result used in the
proof is the Structure Theorem in Kontorovich and Sinai (2002).
Note that the assertion of Theorem 5.3 concerns a
double limit: first $X \to \infty$ and then $m \to \infty$.
See Lagarias and Soundararajan (2006) 
for related results.

\item
Alex V. Kontorovich and Yakov G. Sinai (2002),
{\em Structure Theorem for $(d, g, h)$-maps,}
Bull. Braz. Math. Soc. (N.S.) {\bf 33} (2002), 213--224.
(MR 2003k:11034). \\
\newline
\hspace*{.25in}
This paper studies $(d, g, h)$-maps, in
which $g > d \ge 2$, with $g$ relatively
prime to $d$, 
and $h(n)$ is a periodic
integer-valued function with period $d$, with
$h(n) \equiv -n (\bmod~ d)$ and 
$0 < |h(n)| < g$.
The  $(d, g, h)$-map is defined on $\ZZ \backslash d\ZZ$ by 
$$
T(x) := \frac{gx + h(gx)}{d^k},  ~~\mbox{with}~~ d^k || gx + h(gx).
$$
A path of $m$ iterates can be specified by the values
$(k_1, k_2,..., k_m)$ and a residue class $\epsilon ~(\bmod~dg)$,
and set $k= k_1 + k_2+ ...+ k_m.$ 
The structure theorem states that exactly $(d-1)^m$ triples
$(q, r, \delta)$ with $0 \le q < d^k$, $0 < r < g^m$ and
$\delta \in E= \{ j : 1 \le j < dg, d \nmid j, g \nmid j \}$
produce a given  path, and for such a triple $(q, r, \delta)$
and all $x \in \ZZ$,
$$
T^{(m)} ( gd (d^k x + q) + \epsilon) = gd(g^m x + r) + \delta.
$$ 
They deduce that, as $m \to \infty$, a properly logarithmically scaled version
of iterates converges to a Brownian motion
with drift $\log g - \frac{d}{d-1} \log d$.  More precisely, fix $m$,
and choose points $0= t_0 < t_1 < ...< t_r=1$ and set 
$m_i = \lfloor t_im\rfloor$ and $y_i = \log T^{(m_i)}(x)$. Then the
values $y_i - y_{i-1}$ converge to a Brownian path.
These results imply that when the drift is negative, almost
all trajectories have a finite stopping time with $|T^{(m)}(x)| < |x|.$

\item
Benjamin Kraft and Kennan Monks (2010),
{\em On conjugacies for the $3x+1$ map induced by continuous endomorphisms
of the shift dynamical system,}
Discrete Math. {\bf 310} (2010), no. 13-14, 1875--1883.
(MR 2011e:37027)\\
\newline
\hspace*{.25in} 
The continuous endomorphisms of the shift dynamical system were classified by
Maria Monks (2009). A conjugacy $\Phi$ between  the $3x+1$ map $T$ extended to the $2$-adic
integers to the one-sided shift $S$ was shown in Lagarias (1985), i.e. $\Phi \circ S \circ \Phi^{-1} = T$.
This paper studies maps $H_f = \Phi \circ f \circ \Phi^{-1}$ where $f$ is a continuous
endomorphism of the $2$-adic integers $\ZZ_2$ (a.k. a. binary shift dynamical system).
It generalizes results of Maria Monks (2009).

\item
Ilia Krasikov and J. C. Lagarias (2003),
{\em Bounds for the $3x+1$ problem using difference inequalities},
Acta Arithmetica {\bf 109} (2003), no. 3,  237--258. 
(MR 2004i:11020)\\
\newline
\hspace*{.25in} 
This paper deals with the problem of
 obtaining lower bounds for
$\pi_a(x)$, the counting function for the number of integers
$n \le x$ that have some $3x+1$ iterate $T^{(k)}(n)=a.$ 
It improves the 
nonlinear programming method given 
in Applegate and Lagarias (1995b) for extracting
lower bounds from the inequalities of 
Krasikov (1989). 
It derives a nonlinear program family directly from the
Krasikov inequalities $(\bmod~3^k)$ whose associated lower
bounds are expected to be the best possible derivable by this approach.
The nonlinear program 
for $k=11$ gives the  improved lower bound:
If $a \not\equiv 0~~(\bmod~3)$, then
$\pi_a (x)  >  x^{.841}$
for all sufficiently large $x$. 
The interest of the new nonlinear program family is the 
(not yet realized) hope of proving 
$\pi_a (x)  >  x^{1 - \epsilon}$
by this approach, taking a sufficiently large $k$.

\item
Stuart A. Kurtz and Janos Simon (2007),
{\em The undecidability of the generalized Collatz problem,}
In: J-Y. Cai et al, Ed., {\em Theory and Applications of Models of Computation,
4-th International Conference, TAMC 2007, Shanghai, China, May 22-25, 2007.}
 Lecture Notes in Computer Science No. 4484, Springer-Verlag: New York (2007),
pp. 542--553. (MR 2374341)  \\
\newline
\hspace*{.25in} 
A generalized Collatz
function $g$ is a function mapping the nonnegative 
integers to itself  such that there is a modulus $m$ such that  for
each congruence class $x \equiv i~(\bmod~m)$ it is given by 
$g(x) = a_i x + b_i$ when  where $(a_i, b_i)$ are 
nonnegative rational numbers 
 such that $(i+km) a_i+b_i$ is an integer for all
integer $k$.  Generalized  Collatz functions can
easily  be  listed  with a G\"{o}del numbering $k_e$.

The authors build on the work of Conway (1972), who considered the subclass
of such functions with all $b_i=0$, which we term multiplicative generalized Collatz functions. 
Conway proved that the following problem is
undecidable:  
Given as input  $(g,n)$ where $g$ is  a  multiplicative generalized Collatz function and
$n$  a nonnegative integer, can it be decided  whether there is
 some $i \ge 1$ such that the $i$-th iterate
$g^{(i)}(2^n)$ is a power of $2$?   Conway obtained the result by
showing   that the partially defined function
$M(n)$ determined by such a function $g$ (with $2^{M(n)}$ being the first power of $2$
encountered as an iterate $g^{(i)}(2^n)$ if one exists, and $M(n)$ is undefined otherwise) could
be any unary partial recursive function.
Here the authors formulate:

{\em  Generalized Collatz Problem $(GCP)$.}  Given a representation of a generalized
Collatz function $g$, can it be decided for all $x \ge 0$ whether there exists
some $i \ge 1$ such that $g^{(i)}(x) =1$?

They first formulate (Theorem 1): Let $M_e$ be an acceptable  G\"{o}del
numbering  of partial recursive functions. Then from the index $e$ one can compute
a representation of a generalized Collatz function $g_e$ such that $M_e$ is total if and only
if for every $x \ge 1$ there exist a value  $i ge 1$ such that $g_e^{(i)}(x) =1$.
From this result they formulate (Theorem 2): The problem GCP is $\Pi_0^2$-complete.
This  locates precisely in the degrees of unsolvability hierarchy the difficulty of this
undecidable problem. 

This is a conference paper and does not supply detailed proofs.\\

\item
Jeffrey C. Lagarias (2003+),
{\em The $3x+1$ Problem: An Annotated Bibliography (1963--1999)}, \\
eprint:  {\tt arxiv:math.NT/0309224 Sept. 13, 2003, v11}. \\
\newline
\hspace*{.25in}
This is the initial installment of the annotated bibliography. It contains over
180 items.

\item
Jeffrey C. Lagarias (2006), 
{\em Wild and Wooley Numbers},
American Mathematical Monthly, {\bf 113} (2006), 97--108. (MR 2007g:11029).\\
\newline
\hspace*{.25in}
This paper considers some problems about  multiplicative
semigroups of positive rationals motivated by
work of Farkas (2005) on variants of the 
$3x+1$ problem.  
The {\em wild  semigroup} $\sS$ is the semigroup of
positive rational numbers generated by
$\{ \frac{3n+1}{2n+1}: n \ge 0 \}$ together with $\frac{1}{2}$,
and the {\em Wooley semigroup} is the sub-semigroup
generated by $\{ \frac{3n+1}{2n+1}: n \ge 0 \}$ without $\frac{1}{2}$,
This paper considers the question of which integers occur
in these semigroups. The wild integer
semigroup is the set of all integers in $\sS$, and
generators of the wild integer semigroup are termed
{\em wild numbers}. The paper develops evidence in favor
of the conjecture that the wild numbers consist of the
set of all primes, excluding $3$. It shows that $3$
is not a wild number,  that all other primes
below $50$ are wild numbers, and that there are infinitely
many wild numbers. The term ``wild numbers'' was suggested by  
the novel ``The Wild Numbers'' by Philibert Schogt.
The conjecture above was 
proved subsequently in Applegate and Lagarias (2006).

\item
Jeffrey C. Lagarias and Neil J. A. Sloane (2004),
{\em Approximate squaring},
Experimental Math. {\bf 13} (2004), 113--128.  (MR 2005c:11098).\\
\newline
\hspace*{.25in}
This paper studies iteration of the 
``approximate squaring'' map $f(x) = x \lceil x \rceil$,
and asks the question whether for a rational starting
value $x_0 = r > 1$ some iterate is an integer. It
conjectures that the answer is always ``yes'',  and proves
it for rationals $r$ with denominator $2$. It shows that
this holds for most rationals having a fixed denominator
$d \ge 3$ with an exceptional set of integers below $x$
of size at most 
$O(x^{\alpha_d})$ for certain constants $0 < \alpha_d < 1$.
It then considers a variant of this problem on the
$p$-adic numbers, where an exceptional set exists and
is shown to have Hausdorff dimension equal to $\alpha_p$.

The paper also studies the iteration of
``approximate multiplication'' maps $f_r(x) = r \lceil x \rceil$,
where $r$ is a fixed rational number. It conjectures that for $r >1$ 
all but a finite number of integer starting values have
some subsequent iterate that is an integer, and proves
this for rationals $r$ with denominator $2$. It shows
for rationals $r$ with denominator $d$ that the size of
the exceptional set of integers below $x$ that have
no integer in their forward orbit under $f_r$ has
cardinality at most $O(x^{\beta_d})$ with
$\beta_d = \frac{\log (d-1)}{\log d}.$
It suggests that this conjecture
is likely to be  hard in the general 
case, by noting an analogy with iteration of
the map appearing
in Mahler's $Z$-number problem, see Mahler (1968).

\item
Jeffrey  C. Lagarias and K. Soundararajan (2006),
{\em Benford's Law for the $3x+1$ Function,}
J. London Math. Soc. {\bf 74} (2006), 289--303.
(MR 2007h:37007)\\
\newline
\hspace*{.25in}
Kontorovich and Miller (2005)  proved results 
concerning Benford's law for initial $3x+1$ iterates,
in a double limit as the number of iterates $N \to \infty$.
This paper  proves a quantitative
version of Benford's law valid for finite $N$. 
Benford's law (to base B) for an infinite sequence $\{x_k: k \ge 1\}$
of positive quantities $x_k$ 
is the assertion that $\{ \log_B x_k : k \ge 1\}$ is
uniformly distributed $(\bmod~ 1)$.
This paper studies the initial iterates 
$x_k= T^{(k)}(x_0)$ for $1 \le k \le N$ of the $3x+1$
function, where $N$ is fixed. It shows that for most 
initial values $x_0$, such sequences 
approximately satisfy Benford's law, in the sense that the discrepancy of
the finite sequence $\{ \log_B x_k: 1 \le k \le N \}$ is small.
The precise result treats the uniform distribution of initial
values $1\le x_0 \le X$, with $x \ge 2^N$, and shows that
for any (real) base $B>1$ 
the discrepancy is smaller than  $2 N^{ - \frac{1}{36}}$
for all but an exceptional set $|\sE(X, B)| $ of cardinality
$|\sE(X, B)|  \le c(B) {N^{-\frac{1}{36}}}X$,
where $c(B)$ is independent of $N$ and $X$.

\item
Eero Lehtonen (2008),
{\em Two undecidable versions of Collatz's problems},
Theor. Comp. Sci. {\bf 407} (2008), 596--600. (MR 2009k:68092)\\
\newline
\hspace*{.25in}
The author  gives two constructions. He first constructs  a function  
$$
f(x) = 
\left\{
\begin{array}{cl}
3x+t & \mbox{if}~ n \in A_{t}, ~\mbox{for}~ t= -9, -8, ..., 8, 9, \\
~~~ \\
\df{x}{2} & \mbox{if} ~~x \equiv 0~~ (\bmod~2) ~,
\end{array}
\right.
$$
in which  $\{ A_t: ~-9 \le t \le 9\}$ are recursive sets of nonnegative integers
that partition the class of odd positive integers, for which the problem of
deciding if a given integer iterates to $1$ is undecidable. It encodes the
halting problem for Turing machines,  when halting occurs if and only if the
function  iterations eventually arrive at $1$. His second construction is
a  variant of Collatz's original function
$f(3n)=2n, f(3n+1)=4n+1, f(3n+2)=4n+3$, which is a permutation of
the positive integers $\NN$. He constructs a recursive function 
 $f: \NN \to \NN$ which is a bijection of $\NN$, such that the question
 whether a given input $n$ is in a finite cycle under iteration of $f$ is
 undecidable.

\item
Dan Levy (2004),
{\em Injectivity and Surjectivity of  Collatz Functions},
Discrete Math. {\bf 285} (2004), 190--199. (MR 2005f:11036). \\
\newline
\hspace*{.25in}
This paper gives necessary and sufficient conditions
on members of a class of generalized Collatz maps of the form
$T(x) = \frac{m_i x - r_i}{d}$ for $x \equiv i~(\bmod~d)$
to be injective maps, resp. surjective maps,  on the integers.
These give as a corollary  a criterion of Venturini (1997)
for such a map to be a permutation of the integers.

The author frames some of his results in terms of
concepts involving integer matrices. He introduces
a notion of {\em gcd matrix} if its elements
can be written $M_{ij}=gcd(m_i, m_j)$ and a {\em difference
matrix} if its elements can be written $M_{ij}= m_i - m_j$.
Then he considers a relation that $M$ is a {\em total
non-divisor} of $N$ if $M_{ij} \nmid N_{ij}$ for all
$i,j$. Then the author's condition for injectivity
of a generalized Collatz map above 
is that the $d \times d$ gcd  matrix $M_{ij}= gcd(m_i, m_j)$ 
is a total non-divisor of the $d \times d$ difference matrix
$N_{ij}= q_i- q_j$, with $q_j = \frac{r_j -j m_j}{d}$.

A very interesting result of the author is 
an explicit example of an injective  function $T(\cdot)$
in the class above which
has a (provably) divergent trajectory, and which  has 
iterates both increasing and decreasing in size.
This particular map $T$ is not surjective.

\item
Li, Xiao Chun (2002),
{\em Compressive Iteration for the $3N+1$ conjecture} (Chinese),
J. Huazhong Univ. of Science and Technology (Natural Science)
[Hua zhong gong xue yuan] {\bf 30} (2002), no. 2. \\
\newline
\hspace*{.25in}
[I have not seen this paper.]

\item
Li, Xiao Chun (2003),
{\em Contractible iteration for the $3n+1$ conjecture}
(Chinese, English summary), J. Huazhong Univ. Sci. Technol. Nat. Sci.   
[Hua zhong gong xue yuan] {\bf 31} (2003), no. 7, 
115--116. (MR 2025640).\\
\newline
\hspace*{.25in}
English summary: "The concept of contractible iteration for the $3N+1$ conjecture
was presented. The results of contractible iteration were gien as follows: some
equivalence propositions for the $3N+1$ conjecture; sequence of Syracuse order;
proof of term formula of $n$; some theories about $t_a(n)= t_c(n).$"

{\em Note.} The  adequate stopping time $t_a(n)$ 
and coefficient stopping time $t_c(n)$ are those
given  in  Wu and Hao (2003).
Terras (1976) made a conjecture equivalent to
asserting  the  equality $t_a(n)= t_c(n)$ always holds.

\item
Li, Xiao Chun (2004),
{\em Same-flowing numbers and super contraction iteration in $3N+1$ 
conjecture} (Chinese, English summary),
J. Huazhong Univ. Sci. Technol. Nat. Sci.   [Hua zhong gong xue yuan]  {\bf 32} (2004) no. 10,
30--32.\\
\newline
\hspace*{.25in}
English summary: "The concept of same-flowing numbers in $3N+1$ conjecture
was given. Some definitions and theorems were set up. Research into contraction
iteration of odd numbers $2a+1 (a \in \NN_d$) could be converted to research
into Syracuse order of odd number $a$ and contraction iteration simplified. Infinite
sequence same-flowing to $3a (a \in \NN_d)$ was constructed. The concept
of super contraction iterations and term formula of $x$ was given."

{\em Note.} The super-contraction iteration keeps track of the successive
odd integers in the Collatz iteration.

\item
Li, Xiao Chun (2005),
{\em Necessary condition for periodic numbers in the $3N+1$ conjecture}
(Chinese, English summary), J. Huazhong Univ. Sci. Technol. Nat. Sci.  
[Hua zhong gong xue yuan]  {\bf 33} (2005), no. 11,  102--103. (MR 2209315).\\
\newline
\hspace*{.25in}
English summary: "After the introduction of the conception of periodic numbers
in $3N+1$ conjecture, numerical theory function was defined, the definition of
sum-line number in binary bit described and the Gauss function given.
At the same time, some properties of number theory function and calculation
formula of were discussed. By the use of number theory function, a necessary
condition of periodic numbers in $3N+1$ conjecture was proposed. The presentation
of this necessary condition could play a certain role in further solution
to periodic numbers question in $3N+1$ conjecture."

\item
Li, Xiao Chun  (2006),
{\em Some properties of super contraction iteration in the $3N+1$ conjecture}
(Chinese, English summary), 
J. Huazhong Univ. Sci. Technol. Nat. Sci.  [Hua zhong gong xue yuan] 
{\bf 34} (2006), no. 8, 115--117.
(MR 2287431).\\
\newline
\hspace*{.25in}
English summary: "In researching of $3N+1$ conjecture, by using operator of dividing 
even factor, this paper provided the concept of super contraction and series of iteration
super contraction iteration, which greatly increased iteration velocity in comparison 
with former contraction iteration. At the same time, this paper first provided the concept
of omission numbers and obtained the cyclic relation between contraction iteration
and super contraction iteration, and obtained no unique property about precursor
numbers of monadic-order to odd number in super contraction iteration, and to the odd 
number $y$ of $4k+3$, obtained no unique property about precursor numbers of monadic-order
to odd number in super contraction iteration. At last, this paper first provided necessary
condition of existence of cyclic numbers in super contraction iteration, which can play
active role in cyclic numbers researching of $3N+1$ conjecture. All presented definitions
and theories can simplify in researching world well-known problem $3N+1$ conjecture
in number theory. This paper also provides new method in $3N+1$ conjecture
continuous researching."

{\em Note.} As in his earlier paper  Li (2004), the author studies odd numbers appearing in the
$3x+1$ iteration. 

\item
Li, Xiao Chun and Liang, You Min (2002)
{\em Research on the Syracuse Operator and its properties} (Chinese),
Journal of Air Force Radar Academy (2002), No. 3.\\
\newline
\hspace*{.25in}
[I have not seen this paper.]

\item
Li, Xiao Chun  and Liu, Jun (2006),
{\em Equivalence of the $3N+1$ and $3N+3k$ conjecture and some
related properties} (Chinese, English summary),
J. Huazhong Univ. Sci. Technol. Nat. Sci. [Hua zhong gong xue yuan]  {\bf 34} (2006), no. 7, 120-121.
(MR 2287654)\\
\newline
\hspace*{.25in}
English summary: "$3N+1$ conjecture in number theory was extended to $3N+3^k$ one.
It was pointed out that $3N+3$ conjecture is equivalent to $3N+3^k$ one. Some properties
related to $3N+3$ were obtained. The generalization of $3N+1$ conjecture and these newly
obtained properties not only simplify the operation about $4K+3$-odd numbers in
$3N+3$ conjecture but also provide new way in researching of $3N+1$ conjecture."

\item
Li, Xiao Chun and Wu, Jia Bang  (2004),
{\em Study of periodic numbers in the $3N+1$ conjecture} (Chinese, English summary),
 J. Huazhong Univ. Sci. Technol. Nat. Sci.  [Hua zhong gong xue yuan]  {\bf 32} (2004), 
no. 10, 100--101. (MR 2121229)\\
\newline
\hspace*{.25in}
English summary: "A necessary condition about the existence of periodic numbers
in the $3N+1$ conjecture: $S^{l}(mx_i) =nx_i$, and its generalized necessary
condition: $S^{l} (\sigma_b  x_i)= n^{b+1} x_i$ were presented. A formula about
periodic numbers was given, where $x_1= \frac{r_1}{1- 3^{l}/2^k}$. It was
proved that the length of periodic numbers being $2$ or $3$ could not exist."

{\em Note.} The numbers $x_i$ are the successive odd numbers appearing in
the $3x+1$ iteration. The authors exclude periodic orbits on 
the odd itntegers, where there are exactly
$2$ or $3$ odd numbers. 

\item
Josefine L\'{o}pez and Peter Stoll (2009),
{\em The $3x+1$ conjugacy map over a Sturmian word,}
Integers {\bf 9} (2009), A13, 141--162. (MR 2506145).\\
\newline
\hspace*{.25in}
This paper studies the $3x+1$ conjugacy map $\Phi: \ZZ_2 \to \ZZ_2$
on the $2$-adic integers, which conjugates the $3x+1$ map $T(x)$ to
the $2$-adic shift map, in the sense $\Phi^{-1} \circ T \circ \Phi  = S$,
and which was studied by Bernstein and Lagarias (1996).
It is unknown whether there is any aperiodic $x \in \ZZ_2$ such
that $\Phi(x)$ is periodic; this is conjectured not to happen. 
This paper studies this function for $x$ whose $2$-adic expansion
is a Sturmian word, corresponding to a line with irrational slope. 
This paper finds a generalized continued fraction expansion for 
$\frac{-1}{\Phi(x)}$ in this case convergent in the metric on the $2$-adic integers $\ZZ_2$.
It explicitly computes a number of examples, suggesting that the
images $\Phi(x)$ then have $2$-adic expansions of full complexity.

\item
Florian Luca (2005),
{\em On the nontrivial cycles in Collatz's problem},
SUT  Journal of  Mathematics {\bf 41} (2005), no. 1, 31--41  (MR 2006e:11034). \\
\newline
\hspace*{.25in}
This paper establishes conditions on non-trivial cycles on the positive integers
of the $3x+1$ function. Let $n$ denote the cycle length, and let
$\{x_1,..., x_k\}$ denote the set of odd integers   that appear in such a cycle,
so that  $1 \le k < n$, and let $l_j$ denote the number of
iterates between $x_i$ and $x_{i+1}$, so that  $x_{i+1} = \frac{3x_i+1}{2^{l_i}},$
and $n = \sum_{i=1}^k l_i$.
Let $1 \le J \le k$ denote the number of blocks of consecutive
$l_j$  taking the same value, say $ L_j$, with the block length $N_j$ for
$1 \le j \le J$, so
that $l_i =l_{i+1}=...= l_{i+N_j -1}=L_j$ with $l_{i-1} \ne L_j, l_{i+N_j} \ne L_j$.
Then $L_1 N_1+...+L_j N_j = n$. Call an {\em ascent} a 
value $L_j=1$, i.e. a consecutive string of increasing values 
$x_i < x_{i+1}< ...< x_{i+N_j -1}$ with $x_{i-1}>x_i$ and $x_{i+N_j}< x_{i+N_j-1}$.
The author's main result (Theorem 1) states that there is  an absolute constant $C_1$ such
that there are at least $C_1 \log n$ ascents in any cycle.
In particular for any fixed $c$,
there are only finitely many nontrivial cycles of positive integers having at
most $c$ ascents. This result improves on that of Mimuro (2001).
The proof uses transcendence results coming
from linear forms in logarithms.

This result complements a result of  Brox(2000), who showed
that there are only finitely many integer cycles having at most $2 \log k$
{\em descents}, where a descent is a value $i$ such that $x_{i+1}< x_i$.

{\em Note.} A misprint occurs in the the definition of $l_i$ on
page 32, in condition (ii) "largest" should read "smallest".

\item
Maurice Margenstern (2000),
{\em Frontier between decidability and undecidability: a survey},
Theor. Comput. Sci. {\bf 231} (2000), 217--251.
(MR 2001g:03079). \\
\newline
\hspace*{.25in}
This paper surveys results concerning
the gap between decidability and undecidability,
as measured by the number of states or symbols used in a Turing 
machine,  Diophantine equations, the word problem, molecular computations. 
Post systems, register machines,  neural networks and cellular
automata. The $3x+1$ function
(and related functions) have been studied as a test for possible
undecidability of behavior for small state Turing machines,
which are so small they are not known to simulate a universal
Turing machine. The author summarizes known results in
Figure 12, summarizing how large a Turing machine has to be
(number of states, number of symbols) to encode the $3x+1$
problem. For related results, see Margenstern and Matiyasevich (1999)
and Michel (2004). 

{\em Note.} This paper  is the journal version of a conference paper 
in: M. Margenstern, Ed.,  International Colloquium on Universal Machines and
Computations, MCU '98 , Metz, France, March 23-27,1998, Proceedings,
Volume I, IUT: Metz 1998.

\item
Keith R. Matthews (2005+),
{\em The generalized $3x+1$ mapping,}
 preprint, 23pp., dated Oct. 31, 2005, downloadable 
as pdf file from: 
{\tt http://www.maths.uq.edu.au/$\sim$krm/interests.html} \\
\newline
\hspace*{.25in}
This paper discusses the behavior of $3x+1$-like mappings and 
surveys many example
functions as considered in Matthews and Watts (1984, 1985),
Leigh (1986), Leigh and Matthews (1987) and
Matthews (1992),
see also Venturiri (1992).
A particularly tantalizing example is
$$
U(x) = \left\{
\begin{array}{cl}
7x+3 & \mbox{if}~~ x \equiv 0 ~ (\bmod ~3) \\
~~\\
\df{7x+2}{3} & \mbox{if}~~ x \equiv 1 ~ (\bmod~3) \\
~~~ \\
\df{x-2}{3} & \mbox{if}~~ x \equiv 2 ~ (\bmod~3) ~.
\end{array}
\right.
$$
Almost all trajectories contain an element $n \equiv 0 ~(\bmod~3)$ and once
a trajectory enters the set
$\{n : n \equiv 0 ~(\bmod~3)\}$
it stays there.
Matthews offers \$100 to show that if a trajectory has all
iterates $U^{(k)} (x) \equiv \pm 1 ~(\bmod~3)$ then it
must eventually enter one of the cycles $\{1, -1\}$ or
$\{-2, -4, -2\}$.
The paper also considers some maps on the rings of integers of an
algebraic number field, for example
$U: \ZZ [ \sqrt 2 ] \ra \ZZ [ \sqrt 2 ]$ given by
$$
U( \alpha ) = \left\{
\begin{array}{cl}
\df{(1- \sqrt 2) \alpha}{\sqrt 2} & \mbox{if}~~ 
\alpha \equiv 0 ~(\bmod ~( \sqrt 2 )) ~, \\
~~~\\
\df{3 \alpha+1}{\sqrt 2} & \mbox{if} ~~ 
\alpha \equiv 1 ~ (\bmod ( \sqrt 2 ))~.
\end{array}
\right.
$$
The author conjectures that if $U^{(k)} ( \alpha ) = x_k + y_k \sqrt 2$ is a
divergent trajectory, then
$x_k / y_k \ra - \sqrt 2$ as $k \ra \In$.

\item
Karl Heinz Metzger (2000),
{\em Untersuchungen zum $(3n+1)$-Algorithmus.
Teil II: Die Konstruktion des Zahlenbaums,}
PM (Praxis der Mathematik in der Schule) {\bf 42} (2000), 27--32. \\
\newline
\hspace*{.25in}
The author treats the Collatz function $C(x)$, and studies
the directed  graph ("number-forest") on the positive integers
formed with edges iterating the Collatz function. This graph might
have many connected components. The $3x+1$
Conjecture asserts this graph forms a single tree, with one extra directed edge added 
making the trivial cycle at the bottom.  The author constructs
the graph as follows. To each odd number the author forms a
directed path containing all numbers  $\{2^k u: k \ge 0\}$ (a "Spross"), with directed edges from
$2^{k+1}u \to 2^k u$. He then arranges an infinite  table of all positive integers, in which 
the odd numbers $u$  are placed  successively  in a  bottom
horizontal  row and the "Sprossen"  are then put in vertical
columns over them. To the "Sprossen" edges  he now adds directed edges
starting  from the odd nodes $u$
at the bottom row taking $u \to 3u+1$.
If $u \equiv 1~ (\bmod~4)$ these edges move horizontally to the left and upward
at least two rows, ("leftmovers")  while  if $u \equiv 3~ (\bmod~4)$ they move  horizontally to the right and
upward exactly one row ("rightmovers"). He then considers the successive
odd integers reached during an iteration.
He shows that any  chain of successive
rightmoving odd numbers reached under iteration 
must be finite (Satz 9), i.e. a number of form $4n+1$ is eventually reached
in the iteration. Satz 10 then asserts that positive integers can all be  accounted for  in a 
rearranged  table having only odd numbers of form $4n+1$ along its bottom row.
The results through Satz 9 are rigorous; maybe Satz 10 too.

The final Satz 11 asserts the truth of the $3x+1$ Conjecture.
However its proof  on page 32 is incomplete. The proof seems to implicitly assume the 
 decrease in size during  leftmoving steps overcomes 
 the increase in size of iterates during rightmoving steps, but this is not rigorously shown.

\item
Karl Heinz Metzger (2003),
{\em Untersuchungen zum $(3n+1)$-Algorithmus,
Teil III: \\
Gesetzm\"{a}ssigkeiten der Ablauffolgen,}
PM (Praxis der Mathematik in der Schule) 
{\bf 45} (2003), No. 1, 25--32. \\
\newline
\hspace*{.25in}
This paper is independent of  parts I and II, which had gaps in some proofs, and does not
reference any of the papers Metzger (1995) (1999) (2000).

It studies iterates of the Collatz function $C(n)$, viewing the 
 iterates $(\bmod~6)$. It observes a kind of self-similar structure
 in the tree of inverse iterates, and introduces
 some symbolic dynamics to describe it. It gives formulas for
 elements on certain branches of the tree. If one writes 
 $n=6(k-1)+i$ with $0 \le i \le 5$ and $C(n)= 6(\bar{k}-1)+j$
 then it describes how to update $k$ to $\bar{k}$, over a 
 sequence of iterates.

\item
Pascal Michel (2004),
{\em Small Turing machines and generalized busy beaver competition,}
Theor. Comp. Sci. {\bf 326} (2004), 45--56. [MR 2005e:68049, Zbl 1071.68025]\\
\newline
\hspace*{.25in}
Let $TM(k, l)$ denote the set of one-tape Turing machines with $k$ states and
$l$symbols; the tape is two-way infinite. Much work has been done
on how large  $k$ and $l$ must be to encode an undecidable halting problem.
It is known that the halting problem is decidable if $k=1$ or $l=1$ 
and also for classes $T(2,3), T(3,2)$ (hence $T(2,2)$. 
 It is known that universal Turing machines exist in the classes
$T(2, 18), T(3,9), T(4,6),T(5,5), T(7,4), T(10,3)$ and $T(19,2)$, hence the halting
problem is undecidable in these classes. (See Y.Rogozhin,
Small universal Turing machines, Theor. Comp. Sci. {\bf 168} (1996), 215--240.)
It is also known that the $3x+1$ problem can be
encoded in a Turing machine in classes  $T(2,8), T(3,5), $ $ T(4,4),$ $ T(5,3)$ and $T(2,10)$.
It follows that 
there is no method currently known for solving the halting problem for
Turing machines in these classes. (See Margenstein (2000) for a survey,
and also Michel (1993).)
Here the author shows there are Turing machines computing $3x+1$-like functions
for which halting is not known in the small Turing classes $T(2,4), T(3,3), T(5,2)$.
Thus the  only remaining class for which decidability of the halting problem is likely
to be possible is $T(4,2)$. 

For the class $T(2,4)$ the author constructs a machine that
can encode iteration of  the function 
$g: \ZZ_{\ge 0} \times \{0,1\} \to \ZZ_{\ge 0}  \times \{0,1\} $ 
given by
$$
\begin{array}{lcc}
g(3k, ~~~0) &= &(5k+1, 1)\\
g(3k+1,0)& = &\mbox{halt}\\
g(3k+ 2,0)& = &(5k+4,0) \\
g(3k, ~~~1) &= &\mbox{halt}\\
g(3k+1, 1) &=& (5k+5,0)\\
g(3k+2,1)&=& (5k+7,1).
\end{array}
$$
It is not known whether or not every input value $\ZZ_{\ge 0} \times \{0, 1\}$
iteration of the function $g$ 
eventually halts. 

The paper also constructs machines in small Turing classes
that achieve new records on the maximal number of steps before halting on
empty input (busy beaver function $S(k, l)$) and for  the maximal number
of symbols printed before halting on empty input ($\Sigma(k,l)$. These
show $S(2,3) \ge 38,~S(2,4) \ge 7195,~S(3,3) \ge 40737$
and $\Sigma(2,3) \ge 9, \Sigma(2,4) \ge 90$ and $\Sigma(3,3) \ge 208.$ 
The author conjectures the values for $(2,3)$ and $(2,4)$ are 
best possible.

\item
Tomoaki Mimuro (2001),
{\em On certain simple cycles of the Collatz conjecture,}
SUT Journal of Mathematics, {\bf 37}, No. 2 (2001), 79--89.
(MR 2002j:11018). \\
\newline
\hspace*{.25in}
The paper shows there are only finitely many positive integer cycles 
of the $3x+1$ function whose
symbol sequence has the form  $ 1^i~(10^j )^k $, where $i, j, k$ vary over
nonnegative integers. (The symbol seqence is read left to right.)
This result includes  the trivial cycle starting from
$n=1$, whose  symbol sequence is $(10)$, where $(i, j, k) = (0, 1, 1)$.
Suppose that
 the periodic orbit has period $p=i+k(j+1)$ terms, of which $d= i+k$
are odd. The author shows by elementary arguments that there
are no integer orbits of the above type with 
$\frac{3}{4} > \frac{3^{d}}{2^{p}}$,
and the trivial cycle is the unique solution with 
$\frac{3}{4} =\frac{3^{d}}{2^{p}}$. Using bounds from
transcendence theory (linear forms in logarithms) he
shows that there are finitely many values of  $(i, j, k)$ giving
an integer orbit with $1 >  \frac{3^{d}}{2^{p}} > \frac{3}{4},$
with an effective bound on their size. Any orbit on the
positive integers necessarily has $1 >  \frac{3^{d}}{2^{p}},$
so the result follows.

For other papers using transcendence theory to classify some
types of periodic orbits, see Steiner (1978), Belaga and Mignotte (1999),
Brox (2000), Simons (2005), Simons and de Weger (2005). 
A further improvement on this result is given in Luca (2005).

{\em Note.}
There are two known integer orbits on the negative integers
of the author's  form. They are 
 $n=-1$ with symbol seqence $(1)$, where  $(i, j, k) = (1, *, 0)$,
and  $\frac{3^{d}}{2^{p}}= \frac{3}{2}$, 
and $n=-5$ with symbol sequence  $(110)$ where $(i, j, k) = (1, 1, 1)$,
and $\frac{3^{d}}{2^{p}}= \frac{9}{8}$. 
The author's finiteness result might 
conceivably be extended to  
the range $\frac{3}{2} \ge \frac{3^{d}}{2^{p}} \ge 1$
and so cover them.

\item
Michal Misiurewicz and Ana Rodriguez (2005),
{\em Real $3X+1$},
Proc. Amer. Math. Soc., {\bf 133}  (2005), 1109--1118.  (MR 2005j:37011).\\
\newline
\hspace*{.25in}
The authors consider the semigroup generated by the two maps 
$T_1(x) = \frac{x}{2}$ and $T_2(x) = \frac{3x+1}{2}$.
They show this semigroup is a free semigroup on two
generators.
The forward orbit of a positive  
input $x_0$  under this semigroup is 
$$O^{+}(x_0) 
:= \{ T_{i_1} \circ ...\circ T_{i_n}(x_0) : n \ge 1,~ \mbox{each}~
i_k \in \{0,1\} \}.
$$
They show that each  orbit $O^{+}(x_0)$ is dense on $(0, \infty)$.
Furthermore they show that starting from $x_0$ one can get an
iterate $T^{(n_0)}(x_0)$ within a given error $\epsilon$
of  a given value $y$ while remaining
in the bounded region 
$$
\min(x, y- \epsilon \le T_{i_1} \circ \cdots T_{i_j}(x_0) \le \max (11x+4, 4y -x). ~~ 1 \le j \le n_0.
$$
They show that orbits having a periodic point are dense in $(0, \infty)$.
Finally they show that the group of homeomorphisms of the line
generated by $T_1, T_2$ consists of all
maps 
$x \mapsto 2^k 3^l x + \frac{m}{2^i 3^j} $,
in which $k, l, m$ are integers and $i, j$ are nonnegative. It 
is not a free group.
These results concern topological dynamics;  for results concerning measurable
dynamics see Bergelson, Misiurewicz and Senti (2006).

\item
Kenneth G. Monks (2002), 
{\em $3X+1$ Minus the $+$},
Discrete Math. Theor. Comput. Sci. {\bf 5} (2002), 47--53.
(MR 2203f:11030). \\
\newline
\hspace*{.25in}
This paper formulates a FRACTRAN program (see Conway(1987))
of the form $R_i(n)\equiv r_i n$ if $n \equiv i (\bmod~d)$,
such that the $3x+1$ Conjecture is true if and only if the
$R$-orbit of $2^m$ contains $2$, for all positive integers $m$.
The author determines information on the behavior under iteration of the function
$R(n)$ for all positive integers $n$, not just powers of $2$.
He then deduces information on the possible structure of an integer
$3x+1$ cycle (for the function $T(\cdot)$), 
namely that the sum of its even elements must
equal the sum of its odd elements added to  the number of its odd
elements. Finally he notes that this fact can be deduced directly without 
using the FRACTRAN encoding.

\item
Kenneth G. Monks and Jonathan Yazinski (2004),
{\em The Autoconjugacy of the $3x+1$ function,}
Discrete Mathematics {\bf 275} (2004), No. 1, 219--236. 
MR2026287 (2004m:11030).\\
\newline
\hspace*{.25in}
This paper studies the iteration of the $3x+1$ map $T(x)$
on the $2$-adic integers $\ZZ_2$. It shows that 
the set of 
Aut$(T) = \{ U \in \mbox{Aut}(\ZZ_2): U T U^{-1} = T \}$
consists of the identity map and a map $\Omega= \Phi \circ V \circ \Phi^{-1}$
where $V(x) = -1-x$ is the map reversing the bits in a $2$-adic
integer and $\Phi$ is the $3x+1$ Conjugacy map studied in
Bernstein and Lagarias (1996). It formulates the Autoconjugacy
Conjecture that $\Omega(\QQ_{odd}) \subseteq \QQ_{odd}$,
and proves this  conjecture is equivalent to no rational number with
odd denominator  having a divergent $T$-orbit.
It defines a notion of self-conjugate cycle under the $3x+1$ map,
which is a periodic  orbit $C$ such that $\Omega(C) = C$. It
proves that $\{1, 2\}$ is the only self-conjugate 
cycle of integers. It shows that all self-conjugate cycles
consist of positive rational numbers.

\item
Kenneth M. Monks (2006),
{\em The sufficiency of arithmetic progressions for the $3x+1$ conjecture,}
Proc. Amer. Math. Soc. {\bf 134} (2006), No. 10, 2861--2872.
MR2231609 (2007c:11030).\\
\newline
\hspace*{.25in}
This paper shows that
the $3x+1$ conjecture is true if it is true for all the integers in any
arithmetic progression $\{ A+Bn: n \ge 0\}$, provided $A \ge 0, B \ge 1$.
It gives analogous reductions for the divergent orbits conjecture and
the nontrivial cycles conjecture.

\item
Maria Monks (2009),
{\em Endomorphisms of the shift dynamical systems, discrete derivatives,
and applications,} Discrete Math. {\bf 309} (2009), 5196--5205.
(MR 2010i:37217)\\
\newline
\hspace*{.25in}
The continuous endomorphisms of the one-sided shift dynamical system 
$S$ on the $2$-adic integers are induced by block maps from binary sequences
of length $n$ to those of length $1$. Let $D$ denote the endomorphism
associated to $00 \mapsto 0, 01 \mapsto 1, 10 \mapsto 1, 11 \mapsto 0$.

Also let  $V(x)= 1-x$.  A map $h: \ZZ_2 \to \ZZ_2$ is {\em solenoidal}  if for each $k \ge 1$ 
the first $k$ digits of the input vector $x$ uniquely determine the first $k$ digits of
the output vector $h(x)$.
The main result is  the only continuous endomorphisms of $S$ whose
parity vector is solenoidal are $D, V \circ D, S, $ and $V \circ S$. 

The map  $D$ is interpreted as a ``discrete derivative," and its dynamics
under iteration are studied.  The $2$-adic integers which under iteration
reach a fixed point are characterized. It is shown that each integer $N \in \ZZ$
is eventually periodic under iteration of $D$, with minimal period a power of $2$.

 Finally, a necessary and sufficient condition for the truth of the $3x+1$ conjecture
is given in terms of these concepts.  There is a unique conjugacy $R$  taking $D$ to the
$3x+1$ map $T$, i.e $R \circ D \circ R^{-1} = T$. In Theorem 5.2 the  $3x+1$ conjecture is shown
equivalent to the action of $R^{-1}$ on positive integers $m$ being eventually periodic 
with preperiod of length $2$ and period of length a power of $2$.
\item
Helmut M\"{u}ller (2009),
\:{U}ber Periodenl\"{a}ngen and die Vermutungen von Collatz und Crandall,
[On period lengths and the conjectures of Collatz and Crandall],
Mitt. Math. Ges. Hamburg, {\bf 28} (2009), 121--130.\\
\newline
\hspace*{.25in}
[I have not seen this paper.]
\item
Tomas Oliveira e Silva (2004+),
{\em Computational verification of $3x+1$ conjecture,}
web document at {\tt http://www.ieeta.pt/ \~ tos/}; email: {\tt tos@ieeta.pt}. \\
\newline
\hspace*{.25in}
In Oliveira e Silva (1999) the author reported on computations
verifying the $3x+1$ 
conjecture for $n < 3 \cdot 2^{53} = 2.702 \times 10^{16}$. In 2004 he
implemented an improved version of this
algorithm. As of  February 2008 his computation verified the 3x+1 conjecture 
up to $17 \cdot 2^{58} > 4.899 \times 10^{18}$. 
This is the current record value for verifying the $3x+1$ conjecture.
Compare Roosendaal (2004+).

\item
Reiko Ohira  and Michinori Yamashita (2004),
{\em A Generalization of the Collatz problem} (Japanese), 
PC Literacy [Pasocon Literacy] (Personal Computer Users Application Technology
Association]  {\bf 31} (2006), No. 4, pp. 16--21.\\
\newline
\hspace*{.25in}
The paper presents the $3x=1$ function and then describes
some speed-ups of the iteration, by only stopping at steps
that are odd numbers, with the previous step dividing by
a power $2^e$ with $e>1$. It then defines 
the $p$-Collatz function for an odd $p$ by
$$
f_p(x) = \left\{
\begin{array}{cl}
\df{px + (p-2)}{2}  & \mbox{if} ~~x \equiv 1~~ (\bmod ~ 2) \\
~~~ \\
\df{x}{2} & \mbox{if}~~ x \equiv 0~~ (\bmod ~2 )~.
\end{array}
\right.
$$
The authors study cycles of these functions.
They note  the identity $f_p(p-2)= \frac{(p-2)(p+1)}{2}$. This yields
for $p=2^m-1$ that $f_p(p-2) = (p-2) 2^{m-1},$ so that $p-2$ is
in a periodic orbit; the case $p=3$ gives the trivial cycle of
the $3x=1$ function. It notes that 
 $f_p(p^2-4)= \frac{(p-2)(p-1)^2}{2}.$
This implies 
for  $p=2^m-1$ then $f_p(p^2-4)= (p-2) 2^{2m-1},$  so this enters
the periodic orbit given by $p-2$. 
The paper gives a table of periodic orbits found for various small$p$
including $p= 5, 7,9, 15, 17, 25,27, 29.$ (Note that for $p \ge 5$ most
trajectories of the map $f_p$ are expected to be divergent.)

\item
 Reiko Ohira and Michinori Yamashita(2004),
{\em On the $p$-Collatz problem} (Japanese), 
PC Literacy [Pasocon Literacy] (Personal Computer Users Application Technology
Association] {\bf 31} (2006), No. 4, pp. 61--64.\\
\newline
\hspace*{.25in}
[I have not seen this paper.]

\item
Pan, Hong-liang (2000),
{\em Notes on the $3x+1$ problem} (Chinese),
Journal of Suzhou University, Natural Science 
[Suzhou da xue xue bao. Zi ran ke xue ban]
{\bf 16} (2000), No. 4, 13--16.\\
\newline
\hspace*{.25in}
The author considers the $3x+1$ function $T(n)$,
and first shows the following statement is equivalent to
the $3x+1$ Conjecture: The smallest element $P_n$ 
reached in the  trajectory  of   $n \ge 1$ always
satisfies  $P_n \equiv 1~(\bmod~4)$.

Next, given $n \ge 1$, let $y_k(n)= 1$ (resp. $-1$) according as  
the $k$-th iterate $T^{(k)}(n)$ is odd (resp. even), and define 
$$
D_k(n) := \frac{T^{(k)}(n)}{2^{y_0(n) + y_1(n) + \cdots+ y_{k-1}(n)}}.
$$
The author proves  that  $D_{k+1}(n) \le D_k(n)$ for all $ k\ge 1.$
Since $D_k(n)$ is nonnegative it  follows that the limit $D_n := \lim_{k \to \infty} D_k(n)$ exists.
For initial values  that enter the trivial cycle, one has $D_n >0$, since
eventually the $y_k(n)$ alternate between $+1$ and $-1$. The author 
 shows that  $D_n=0$ for 
all initial values $n \ge 1$  that  either enter a non-trivial cycle or else have  a divergent
trajectory.  It follows that  the $3x+1$ conjecture is equivalent to the assertion 
that $D_n>0$ for all $n \ge 1$.

\item
Joseph L. Pe (2004), 
{\em The $3x+1$ Fractal, }
Computers \& Graphics  {\bf 28} (2004), 431--435. \\
\newline
\hspace*{.25in}
This paper considers iteration of the following
extension of the Collatz function to complex
numbers $z$, which he terms the complex
Collatz function.  Define $C(z) = \frac{z}{2}$
if $\lceil |z| \rceil$ is an even integer, and
$C(z)= 3z+1$ otherwise. A complex number has the
{\em tri-convergence property} if its iterates
contain three subsequences which converge to
$1$, $4$ anbd $2$, respectively. 
The $3x+1$ conjecture  now asserts that
all positive integers have the tri-convergence property.
He gives a sufficient condition for a complex number
to have this property, and uses it to show that
$z= 1+i$  has the tri-convergence property. 
He states that it is unlikely that $z=3+5i$ has this
property. The $3x+1$ problem now asserts that
all positive integers have the tri-convergence property.
He gives some density plots of iterates exhibiting
where they are large or small; self-similarity
patterns are evident in some of them. The author
makes conjectures about some of these patterns,
close to the negative real axis. 

\item
Yuval  Peres, K\'{a}roly Simon and Boris Solomak (2006), 
{\em Absolute continuity for random iterated function systems with
overlaps,}
J. London Math. Soc. {\bf 74}, No. 2 (2006), 739--756.
MR2286443 (2007m:37053). \\
\newline
\hspace*{.25in}
This paper contains results which apply to a question of Ya. G.
Sinal  which was motivated by  the $3x+1$ iteration. It is
stated  in  Section 3, as follows. Let $0<a<1$ be constant 
let   $Z_i$  $(i \ge 1)$ be independent identically distributed discrete
random variables taking  values $1+a$ or
$1-a$, each with probability $\frac{1}{2}.$ Form the random
variable
$$
X := 1+ Z_1 + Z_1 Z_2 +... + Z_1Z_2 \cdots Z_n+ ...
$$
With probability one this sum converges, and the distribution
of $X$ is given by a measure $\nu^{a}$  supported on the interval
$I_a= [ \frac{1}{a}, +\infty)$.  Sinai asked for which values of $a$  is the
measure $\nu^a$ absolutely continuous with 
respect to Lebesgue measure. It is known that when  $\frac{\sqrt{3}}{2} < a < 1$,
the measure $\nu^{a}$  is singular with respect to Lebesgue measure. 
Conjecturally $\nu^a$ is absolutely continuous for almost all $a$ in the
interval $(0, \frac{\sqrt{3}}{2}),$ and this is unsolved.\\

The results of this paper imply the truth of a 
"randomly perturbed" version of this conjecture.  
In Corollary 3.1 the authors consider
$$
X= 1+ Z_1^{'} + Z_1 ^{'}Z_2^{'} +... + Z_1^{'}Z_2^{'}\cdots Z_n^{'}+ ...
$$
in which  $Z_i^{'}:= Z_iY_i$, in which each $Y_i$ is drawn from a fixed absolutely
continuous distribution $Y$ supported on $(1- \epsilon_1, 1+ \epsilon_2)$ 
for small positive $\epsilon_1, \epsilon_2$, having a bounded density 
and  expectation $\EE[ \log Y]=0$. The $Y_i$ and $Z_i$ are drawn independently  
of all other $Y_i$ and all other $Z_i$. Let $\nu_{\bf y}^a$ be the conditional distribution
of $Z_i^{'}$ given a fixed sequence of draws $\by= (y_1, y_2, y_3, ...)$ for all the $Y_i$. 
The conclusions are:

(a) If $0< a< \frac{ \sqrt{3}}{2}$,
then almost every choice of draws $\by=(y_1, y_2, y_3, ....)$ yields an
absolutely continuous conditional distribution $\nu_{\by}^a$. 

(b) If $ \frac{\sqrt{3}}{2} \le a < 1$, then $\nu_{\by}^a$ is always singular with respect to Lebesgue
measure.  For almost all $\by$ the Hausdorff dimension of the support of $\nu_{\by}^a$
is $(2 \log 2)/ \log(\frac{1}{1-a^2}).$

The paper contains many other results.

\item
Qu, Jing-hua (2002),
{\em The $3x+1$ problem and its necessary and sufficient condition} (Chinese),
Journal of Sangqiu Teacher's College (2002), No. 5.\\
\newline
\hspace*{.25in}
[I have not seen this paper.]

\item
Eric Roosendaal (2004+),
{\em On the $3x+1$ problem}, web document, available at: \\
{\tt http://www.ericr.nl/wondrous/index.html} \\
\newline
\hspace*{.25in}
The author maintains an ongoing  distributed search program
for verifying the $3x+1$ Conjecture to new records and for
searching for extremal values for various quantities associated
to the $3x+1$ function. These include quantities
termed the glide, delay, residue,
completeness,  and gamma. Many
people are contributing time on their computers to this
project.

As of February 2008   the $3x+1$ Conjecture is verified up
to $612 \times 2^{50} \approx  6.89 \times 10^{17}.$
The largest value $\gamma(n)$ found so far is
$36.716918$ at 
$n= 7,219,136,416,377,236,271,195 \approx 7.2 \times 10^{21}$.

[The current record for verification of
the $3x+1$ conjecture published in archival literature is that
of Oliveira e Silva (1999). Note that Oliveira e Silva has
extended his computations to $4.899 \times 10^{18},$ the current
record.]

\item
Jean-Louis Rouet and Marc R. Feix (2002),
{\em A generalization of the Collatz problem. Building
cycles and a stochastic approach,}
J. Stat. Phys. {\bf 107}, No. 5/6 (2002), 1283--1298. \\
(MR 2003i:11035). \\
\newline
\hspace*{.25in}
The paper studies the class of functions
$U(x) = (l_i x + m_i)/n$ if $x \equiv i~(\bmod~n)$,
with $i l_i + m_i \equiv 0~(\bmod~n).$ 
These functions include the
$3x+1$ function as a special case. 
They show that there is a bijection between the
symbolic dynamics of the first $k$ iterations
and the last $k$ digits of the input $x$ 
written in base $n$ if and only if $n$ is relatively prime
to the product of the $l_i$.
They show that for fixed $n$ 
and any given $\{ m_i: 1 \le j \le k\}$
and can find a set of coefficients $\{l_i: 0 \le i \le n-1\}$
and $\{m_i: 0 \le i \le n-1\}$
with $n$ relatively prime to the product of the
$l_i$ which give these values as a $k$-cycle,
$U(m_i) = m_{i+1}$ and $U(m_k)= m_0.$
They give numerical experiments indicating that
for maps of this kind on $k$ digit inputs (written in base 
$n$) ``stochasticity'' persists
beyond the first $k$ iterations. 
For the $3x+1$ problem itself (with base $n=2$), it is believed that 
 ``stochasticity''
persists for about $c_0~k$ iterations, with
$c_0 = \frac{2}{\log_2(4/3)} = 4.8187$, as
described in Borovkov and Pfeifer (2000), who also present
supporting numerical data.

\item
Giuseppe Scollo (2005),
{\it $\omega$-rewriting the Collatz problem,}
Fundamenta Informaticae {\bf 64} (2005), 401--412.
[MR 2008g:68060, Zbl 1102.68051] \\
\newline
\hspace*{.25in}
This paper reformulates the Collatz iteration dynamics
as a term rewriting system.
$\omega$-rewriting allows
infinite input sequences and infinite rewriting. 
In effect the dynamics is extended to a larger domain,
allowing infinitary inputs. The author shows the inputs 
extend to the 
the $3x+1$ problem on rationals with odd denominator.
The infinitary  extension seems analogous to extending the $3x+1$ map
to the $2$-adic integers. 

\item
Douglas J. Shaw (2006),
{\em The pure numbers generated by the Collatz sequence},
Fibonacci Quarterly, {\bf 44} , No. 3, (2006), 194-201.\\
\newline
\hspace*{.25in}
A positive  integer is called {\em pure} if its entire tree of preimages under
the Collatz map $C(x)$ contains no integer that is smaller than it is; otherwise
it is called {\em impure}. Equivalently, an integer $n$ is {\em impure} if
there is some $r<n$ with $C^{(k)}(r) =n$ for some $k \ge 1$. Thus $n=4$
is impure since $C^{(5)}(3)=4$. This paper develops congruence  conditions
characterizing pure and impure numbers, e.g. all $n \equiv 0 ~(\bmod~18)$
are impure, while all $n \equiv 9~ (\bmod ~18)$ are pure. 
 It proves  that the set of pure
numbers and impure numbers each have a natural density. The density
of impure numbers satisfies $\frac{91}{162} < d < \frac{2}{3}.$
The subtlety in the structure of the set of pure numbers concerns which numbers
$n \not\equiv 0~(\bmod~3)$ are pure.

\item
Shi, Qian Li (2005a),
{\em Properties of cycle sets} (Chinese),
J. Yangzte Univ. Natural Science.
[Changjiang daxue xuebao. Zi ke ban]
{\bf 2} (2005), No. 7, 191--192. (MR 2241924)\\
\newline
\hspace*{.25in}
[I have not seen this paper.]

\item
Shi, Qian Li (2005b),
{\em Properties of the $3x+1$ problem} (Chinese),
J. Yangzte Univ. Natural Science.
[Changjiang daxue xuebao. Zi ke ban]
{\bf 2} (2005),  No. 10, 287--289. (MR 2241263)\\
\newline
\hspace*{.25in}
[I have not seen this paper.]

\item
John L. Simons (2005),
{\em On the non-existence of $2$-cycles for the $3x+1$ problem},
Math. Comp. {\bf 74} (2005), 1565--1572. MR2137019(2005k:11050).\\
\newline
\hspace*{.25in}
The author defines an {\em $m$-cycle} of the $3x+1$ problem
to be an periodic orbit of the $3x+1$ function
that contains $m$ local minima, i.e. contains
$m$ blocks of consecutive odd integers.  This
is equal to the number of {\em descents} in the  terminology
of Brox (2000), which correspond to  local maxima. 
The result of Steiner (1978) shows there are no nontrivial
$1$-cycles of the $3x+1$ problem on the positive integers. 
The author  proves there are no non-trivial $2$-cycles for
the $3x+1$ function. See Simons and de Weger (2005) for
an impressive generalization  of this result. 

\item
John L. Simons (2007),
{\em A simple (inductive) proof for  the non-existence of $2$-cycles for the $3x+1$ problem},
J. Number Theory  {\bf 123} (2007), 10-17. \\
\newline
\hspace*{.25in}
This paper gives complements  Simons (2005) by giving another
proof of the non-existence of $2$-cycles (periodic orbits containing exactly
two blocks of consecutive odd elements) for the $3x+1$ function on the positive integers.
The last section sketches a proof that the $3x-1$ function on the
positive integers has a single $2$-cycle with minimal element
$n_0=17$; this is equivalent to showing that the $3x+1$ function on
the negative integers has a single $2$-cycle, starting from $n_0=-17$.
The author's method applies to the $3x\pm q$ problem, with fixed $q$ with
$gcd(6,q)=1$ to find a finite list of $2$-cycles; however it does not extend to 
classify $m$-cycles with $m \ge 3$.
For results on $m$-cycles, see  Simons and de Weger(2005).

\item
John L. Simons (2008),
{\em On the (non)-existence of $m$-cycles for generalized
Syracuse sequences, }
Acta Arithmetica {\bf 131} (2008), No. 3, 217--254. (MR 2008m:11056)\\
\newline
\hspace*{.25in}
The author defines an {\em $m$-cycle} of the $3x+1$ problem
to be a periodic orbit of the $3x+1$ function
that contains $m$ local minima, i.e. contains
$m$ blocks of consecutive odd integers. 
This paper genralizes a proof of Simons and de Weger (2005) to give criteria for  the
non-existence of $m$-cycles for the $3x+1$ problem of generalized Collatz
sequences such as the $3x+q$ problem, the $px+1$ problem and the
inverse Collatz problem. (Note that such problems may have $m$-cycles for
various small values of $m$.)

\item
John L. Simons (2008),
{\em Exotic Collatz Cycles,}
Acta Arithmetica  {\bf 134} (2008), 201--209. (MR 2438845)\\
\newline
\hspace*{.25in}
The author considers the $px+q$ problem, with $\gcd(p. q)=1$, $p \ge 5$ odd. 
He studies primitive $m$-cycles. He shows, modulo a conjecture, that for
each $p$ there exist a $q$ with the number of cycles being arbitrarily large.
He calls such examples ``exotic" because empirically such problems appear to 
have most orbits diverging to infinity.

\item
John L. Simons (2008a+),
{\em Post-transcendence conditions for the existence of $m$-cycles for
the $3x+1$ problem,} preprint.\\
\newline
\hspace*{.25in}
The author defines an {\em $m$-cycle} of the $3x+1$ problem
to be a periodic orbit of the $3x+1$ function
that contains $m$ local minima, i.e. contains
$m$ blocks of consecutive odd integers. 
This paper extends the bounds of Simons and de Weger for non-existence
of $m$-cycles for the Collatz problem to larger values of $m$, including
all $m \le 73$.

\item
John L. Simons (2008b+),
{\em On isomorphism between Farkas sequences and Collatz sequences,}
preprint.\\
\newline
\hspace*{.25in}
The author considers generalizations of  sequences
studied by  Farkas (2005). These are sequences taking odd integers
to odd integers, of two types, $F_1(a,b,c,d(x)$ 
and $F_2(a,b,c,d)(x)$ where $a,b,c,d$ are odd integers. They
are given by

$$
F_1(x)=\begin{cases} 
            \frac{ax+b}{2}   & \mbox{if}~x \equiv 1 (\bmod ~4); \\
            \frac{cx+d}{2} & \mbox{if}~x \equiv 3 (\bmod~4).
                      \end{cases}
$$
and
$$
F_2(x)=\begin{cases}
               \frac{x}{3} & \mbox{if}~ x\equiv 3, 9~(\bmod~12);\\
              \frac{ax+b}{2}&  \mbox{if}~ x\equiv 1,5~(\bmod~12);   \\
            \frac{cx+d}{2} & \mbox{if}~ x\equiv 7,11~(\bmod~12).
             \end{cases}
$$
He introduces a notion of {\em isomorphism} between orbits 
$\{ x_n\}$ and $\{y_n\}$ of two recurrence sequences, namely that
there are nonzero integers $\alpha, \beta$ and an integer $\gamma$
such that
$$
\alpha x_k + \beta y_n = \gamma,   ~~\mbox{for~all}~~n \ge n_0.
$$
This notion takes place on the orbit level
and is weaker than that of conjugacy of the maps.
For example, he  observes that each $3x+1$ orbit is  isomorphic to some $3x+q$ orbit.

The author shows that the some sequences $F_1(a,b,c,d)(x)$ 
are isomorphic  to some $px+q$ orbits.  He shows that no 
$F_2(a,b,c,d(x)$ orbit is isomorphic to a $px+1$-map orbit.

\item
John L. Simons and Benne M. M. de Weger (2004),
{\em Mersenne en het Syracuseprobleem}
[Mersenne and the Syracuse problem] (Dutch),
Nieuw Arch. Wiskd. {\bf 5} (2004), no. 3, 218--220.
(MR2090398).\\ 
\newline
\hspace*{.25in}
This paper is a brief survey of the work of
the authors on cycles for the $3x+1$ problem,
given in Simons(2005) and Simons and de Weger (2005).
It also  considers cycles for the $3x+q$ problem, with
$q >0$ an odd number. It defines  
the invariant $S(q)$ to be the sums of the lengths
of the cycles of the $3x+1$ function on the positive
integers. 
The invariant $S(q)$ is not known to be finite
for even a single value of $q$, though the
$3x+1$ conjecture implies that $S(1)=2$.
The paper observes  that  $S(3^t) =S(1)$, for
each $t \ge 1$.
It also considers the case
that $q= M_k:=2^k -1$ is a Mersenne number. 
It gives computational  evidence suggesting that the 
minimum of $S(M_k^2) - S(M_k)$, taken over all $k \ge 3$,
occurs for $k=3$,
with  $S(M_3)= 6$ and $S(M_3^2) = 44$.

\item
John L. Simons and Benne M. M. de Weger (2005),
{\em Theoretical and computational bounds for $m$-cycles of the
$3n+1$ problem}, Acta Arithmetica, {\bf 117} (2005), 51--70.  
(MR 2005h:11049).\\
\newline
\hspace*{.25in}
The authors define an {\em $m$-cycle} of the $3x+1$ problem
to be an orbit of the $3x+1$ function
 that contains $m$ local minima, i.e. contains
$m$ blocks of cosecutive odd integers.  This
is equal to the number of descents in the  terminology
of Brox (2000), which correspond to  local maxima. 
These methods rule out infinite classes of possible symbol
sequences for cycles. 
 The author's main result is  that for each fixed $m$
there are only finitely many $m$-cycles,
there are no 
non-trivial such cycles on the positive integers for
$1 \le m \le 68$, and strong constraints are put on
any such cycle of length at most $72$. The finiteness
result on $m$-cycles was established earlier by Brox (2000).
To obtain their sharp computational
results they  use a transcendence result of
G. Rhin [Progress in Math., Vol 71 (1987), pp. 155-164]
 as well as other methods,
and extensive computations.
[Subsequent to publication, the authors showed there
are also no notrivial such cycles for $69 \le m \le 74$,
see {\tt http://www.win.tue.nl/~bdeweger/3n+1{$\_$}v1.41.pdf}]

\item
Yakov G. Sinai (2003a), 
{\em Statistical $(3X+1)$-Problem},
Dedicated to the memory of J\"{u}rgen K. Moser.
Comm. Pure Appl. Math. {\bf 56} No. 7 (2003), 1016--1028.
(MR 2004d:37007). \\
\newline
\hspace*{.25in}
This paper analyzes iterations of the variant of the
$3x+1$ map that removes all powers of $2$ at each step,
so takes odd integers to odd integers; the author restricts
the iteration to the set $\Pi$ of positive
integers congruent to $\pm 1 (\bmod~6)$,
which is closed under the iteration. It gives a
Structure Theorem for the form of the iterates having
a given symbolic dynamics. The discussion in 
section 5  can be roughly stated as asserting: There is
an absolute constant $c > 0$ such that ``most''  $3x+1$ trees
$(\bmod~3^m)$ contain  at most  $ e^{c \sqrt{m \log m}} 2^{2m}/3^m $ 
nodes whose path to the root
node has length at most $2m$ and which has exactly  $m$ odd iterates.
Here one puts a probability density on such trees which
for a given tree counts the number of such nodes divided by the total
number of such nodes summed over all trees,
and ``most'' means that
the set of such trees having the property contains
$1 - O (1/m)$ of the total probability, as $m \to \infty$.
(The total number of such nodes  is $2^{2m}$, and the total number of trees 
is $2\cdot 3^{m-1}$.) Furthermore at least $\frac{1}{m}$ of
the probability is distributed among trees having at most
$ M^c 2^{2m}/3^m $ such nodes. From this latter result follows
an entropy inequality (Theorem 5.1) which is the author's
main result: The entropy $H_m$ of this probability
distribution satisfies $H_m \ge m \log 3 - O(\log ~m).$
For comparison the uniform distribution on $[1, 3^m]$
has the maximal possible
entropy $H =m \log 3.$
He conjectures that the entropy satisfies  $H_m \ge m \log 3 - O(1).$

See Kontorovich and Sinai (2002) for related results on
the paths of iterates of $3x+1$-like maps.

\item
Yakov G. Sinai (2003b), 
{\em Uniform distribution in the $(3x+1)$ problem},
Moscow Math. Journal {\bf 3} (2003), No. 4, 1429--1440.
(S. P. Novikov 65-th birthday issue). MR2058805 (2005a:11026).\\
\newline
\hspace*{.25in}
Define the map $U(x)$ taking the set $\Pi$ of postive integers congruent to
$\pm 1~(\bmod~6)$ into itself, given by 
$U(x) = \frac{3x+1}{2^k}$ where $k=k(x)$ is the largest power
of $2$ dividing $3x+1$. Then consider all the preimages 
at depth $m$ under $U(\cdot)$ 
of a given  integer 
$y= 6 r + \delta$ with integer $ 0 \le r < 3^m$ and $\delta= \pm 1$.
This consists of 
the (infinite) set of all integers $x$ such that $U^{(m)}(x) =y$.
Let such  a preimage $x$ have  associated
data $(k_1, k_2, ..., k_m)$
with $k_j = k( U^{(j-1)}(x))$, and assign to $x$
the weight $2^{-(k_1+ ... + k_{m})}$ multiplied by $\frac{1}{3}$
if $\delta=1$ and by $\frac{2}{3}$ if $\delta= -1$.
Define the  mass
assigned to $y$ to be  the sum of the weights of all its
preimages at depth m. The sum of
these masses  over
$0 \le r < 3^m$ and $\delta = \pm 1$ 
adds up to $1$ by the Structure Theorem proved
in Sinai (2003).
Let the  scaled size of $y$ be the rational number $\rho = \frac{y}{3^m}$.
This now defines a probability distribution $P(m)$ on these values $\rho$
viewed as a discrete set
inside $[0, 1]$. The main theorem states that as $m \to \infty$
these probability distributions $P(m)$ weakly converge to the uniform
distribution  on $[0,1]$.

\item
Yakov G. Sinai (2004),
{\em A theorem about uniform distribution,}
Commun. Math. Phys. {\bf 252} (2004), 581--588. (F. Dyson birthday issue)
MR2104890 (2005g:37009).\\
\newline
\hspace*{.25in}
This paper presents  a simplified and 
stronger version of the uniform distribution
theorem given in Sinai (2003b).

\item
Matti K. Sinisalo (2003+),
{\em On the minimal cycle lengths of the Collatz sequences},
preprint., Univ. of Oulu, Finland.\\
\newline
\hspace*{.25in}
This paper shows that the minimal length of a nontrivial
cycle of the $(3x+1)$-function on the positive integers
is at least 630,138,897. It uses a method similar to
that of Eliahou (1993), and takes advantage of the 
verification of the $3x+1$ conjecture below the 
bound $ 2.70 \times 10^{16}$ of Oliveira e Silva (1999).
It also considers bounds for cycles of the $(3x -1)$-function.

\item
Alain Slakmon and Luc Macot (2006),
{\em On the almost convergence of Syracuse sequences},
Statistics and Probability Letters {\bf 76}, No. 15 (2006), 1625--1630.\\
\newline
\hspace*{.25in}
The paper  shows that the
 "random "Syracuse conjecture is true in the
sense that random Syracuse sequences get smaller than some
specified bound $B \ge 1$  almost
surely. Consider  identical independent $0-1$ random variables $X_n$ 
having  probability $P[X_n=1]= p$, $P[X_n=0]= q=1-p$.
The authors consider the random Syracuse model $S_{n+1} = \frac{3}{2} S_n + \frac{1}{2}$,
if $X_{n+1}=1$  and $S_{n+1} = \frac{1}{2} S_n$ if $X_{n+1}=0$,
starting from a given $S_0$. 
For the actual $3x+1$ problem one would take $p=q= \frac{1}{2}$. 
They then consider an  auxiliary sequence $Y_n$ with $Y_0= S_0$, formed by the
rule $Y_{n+1} = Y_n \left( 1 + \sigma \gamma \right) $ if $X_{n+1}=1$
and $Y_{n+1} = Y_n \left( 1 - \sigma \right)$ if $X_n=0$, where $\sigma$
and $\gamma$ are positive constants. They formulate results
showing that if $p\gamma- q>0$ then there is a positive threshold value $c$ such
that $Y_n \to \infty$ almost surely if $0< \sigma < c$ and $Y_n \to 0$ almost
surely if $\sigma > c$. 
Finally they show that the 
parameters $(\sigma, \gamma)$ can be chosen so that $\sigma > c$ such
that for 
all  starting values $S_0\ge B$, one has $T_n \ge S_n$ holding at every
step, as long as $T_n \ge B$, while $T_n \to 0$ with probability one. 
The conclude that some $S_n \le B$ with probability one.

\item
Bart Snapp and Matt Tracy (2008),
{\em The Collatz problem and Analogues,}
J. Integer Sequences {\bf 11} (2008), Article 08.4.7, 10pp. (MR 2009i:11144)\\
\newline
\hspace*{.25in}
The authors study a polynomial analogue of the Collatz problem,
as considered in Hicks et al (2008), and then give an
application to the original $3x+1$ problem on $\ZZ$.  They consider for polynomials
in $\ZZ/n\ZZ[X]$ the map
$$
T_p(f) = \left\{
\begin{array}{cl}
\df{(x+1)f(x)- f(0)}{x}  & \mbox{if} ~ f(x) ~\mbox{is not divisible by}~ x \\
~~~ \\
\df{f(x)}{x} & \mbox{if}~~f(x) ~\mbox{is divisible by}~ x.
\end{array}
\right.
$$
They show that for $n=2$ that all polynomials eventually iterate to $1$ under this
map, but for $n \ge 3$ there always exists some polynomial $f(x)$ that never iterates
to $1$ under this map. In Section 3 they apply this result for $n=2$ to define an
enlarged version of the Collatz graph on (all) integers which they call the
{\em combined Collatz graph}. This graph glues together positive
and negative integer iterates, adding solid and  dotted edges, representing
moves of  the $3x+1$ and  $3x-1$ iterations., respectively. They observe that  the $3x+1$
conjecture is equivalent to the claim that all positive integers appear in this graph
connected by undotted edges.

\item
Tang, Ren Xian (2006),
{\em About recursion relation of $3x+1$ Conjecture (Chinese),}
J. of Hunan Univ. Science Engineering 
[Hunan ki ji da xue xue bao] 
(2006), No. 5 \\
\newline
\hspace*{.25in}
[I have not seen this paper.]

\item
Toshio Urata (2002),
{\em Some holomorphic functions connected with the Collatz problem,}
Bulletin Aichi Univ. Education (Natural Science)
[Aichi Kyoiku Daigaku Kenkyu hokoku. Shizen kagaku.] {\bf 51}  (2002), 13--16.\\
\newline
\hspace*{.25in}
The author constructs an entire function $F(z)$ which interpolates the
values of the speeded up Collatz function $\phi(n)$ which takes odd integers to odd
integers by dividing out all powers of $2$, i.e. for an odd integer $n$, 
$\phi(n) = \frac{3n+1}{2^e(3n+1)}$,
where $e(m)= ord_2(m).$  Furthermore the entire  function $F(z)$ is constructed
so that $F'(z)=0$ at all positive odd integers. The author analyzes the holomorphic
dynamics of iterating $F(z)$. He observes that $z=1$ is a superattacting fixed
point, that $z=0$ is a repelling fixed point, and that there is an attracting fixed
point $z$ on the negative real axis with $-\frac{1}{20} < z<0$. He proves that all positive
odd integers are in the Fatou set of $F(z)$, and that all negative odd integers are in the
Julia set of $F(z)$. He observes that the intersection of the Julia set with the negative
real axis coincides with the inverse iterates of $z=0$. Finally he notes that
every component of the Fatou set is simply connected. 

\item
Toshio Urata (2003),
{\em  The Collatz problem over $2$-adic integers,}
Bulletin Aichi Univ. Education (Natural Science)
[Aichi Kyoiku Daigaku Kenkyu hokoku. Shizen kagaku.] {\bf 52}  (2003), 5--11.\\
\newline
\hspace*{.25in}
This is an English version of Urata (1999).
The author studies a $2$-adic interpolation of the speeded-up Collatz function $\phi(n)$
defined on odd integers $n$ by 
dividing out all powers of $2$, i.e. for an odd integer $n$, 
$\phi(n) = \frac{3n+1}{2^p(3n+1)}$,
where $p(m)= ord_2(m).$  Let $\ZZ_2^{\ast}=\{ x \in \ZZ_2: ~x \equiv 1 ~(\bmod~2)\}$ denote the
$2$-adic units. The author sets $OQ:= \QQ \cap \ZZ_2^{\ast}, $ and one has
$\ZZ_2^{\ast}$ is the closure $\overline{OQ}$ of $OQ \subset \ZZ_2$. 
The author shows that the map $\phi$ uniquely
extends to a continuous function
 $\phi: \ZZ_2^{\ast} \backslash \{\frac{-1}{3} \}  \to \ZZ_2^{\ast}$.
He shows that if $f(x) = 2x+ \frac{1}{3}$ then $f(x)$ leaves $\phi$ invariant,
in the sense that $\phi(f(x))= \phi(x)$ for all $x \in \ZZ_2^{\ast} \backslash \{\frac{-1}{3} \} .$
It follows that $f(f(x))= 4x+1$ also leaves $\phi$ invariant.

To each $x \in \ZZ_2^{\ast}  \backslash \{\frac{-1}{3} \} $ he associates
 the sequence of $2$-exponents  $(p_1, p_2, ...)$ produced by iterating $\phi$.
 He proves that
an element $x \in \ZZ_2^{\ast}  \backslash \{\frac{-1}{3} \} $ uniquely determine $x$;
and  that every possible sequence corresponds to some value
$x \in \ZZ_2^{\ast}  \backslash \{\frac{-1}{3} \} $ 
He shows that all periodic points of $\phi$ on $\ZZ_2^{\ast}$ are rational numbers
 $x=\frac{p}{q} \in OQ$,and that there
 is a unique such periodic point for any finite sequence $(p_1, p_2, \cdots, p_m)$
 of positive integers, representing $2$-exponents, having  period $m$.
 If $ C(p_1, p_2, ..., p_{m}) = \sum_{j=0}^{m-1} 2^{p_1 + \cdots + p_j} 3^{m-1-j} $
 then this periodic pont is 
 $$
x= R(p_1, p_2..., p_m):= \frac{C(p_1, ..., p_m)}{2^{p_1 + \cdots + p_m} - 3^m}
$$
He shows that an orbit is periodic if and only if its sequence of $2$-exponents is
 periodic. Examples are given.


\item
Toshio Urata (2007),
{\em Collatz's problem} (Japanese), 
Aichi University of Education Booklet, Math. Sciences Select 1, 
60 pages.\\
\newline
\hspace*{.25in}
This short book considers many aspects of the $3x+1$
iteration, as discussed in his papers
Urata (2002) (2003), (2005). In particular,
it deals with holomorphic dynamics of iterating a function
that is a speeded-up version of the Collatz function.
It includes many pictures of parts of the Fatou and Julia set
for these dynamics.

\item
Toshio Urata and Kazuhiro Hamada (2005),
{\em  Positive values of a holomorphic function connected with the Collatz problem,}
Bulletin Aichi Univ. Education (Natural Science), 
[Aichi Kyoiku Daigaku Kenkyu hokoku. Shizen kagaku.] {\bf 54}  (2005), 1--10\\
\newline
\hspace*{.25in}
The authors study the entire function $F(z)$ introduced
in Urata (2002) which interpolates 
the speeded-up Collatz function $\phi(n)$ defined on the odd integers by 
$\phi(n) = \frac{3n+1}{2^{p(3n+1)}}$,
where $p(m):= ord_2(m).$ They show that $F(x)>0$ for all real $x>0$.
This holds even though $F(z)$
wildly oscillates on the positive real axis.

\item
Tanguy Urvoy (2000),
{\em Regularity of congruential graphs},
in: {\em Mathematical Foundations of Computer Science 2000
(Bratislava)}, Lecture Notes in Computer Science Vol. 1893,
Springer: Berlin, 2000, pp. 680--689. 
(MR 2002d:68083). \\
\newline
\hspace*{.25in}
The Collatz problem is viewed as a reachability problem on
a directed graph with vertices labelled by positive 
integers $n$, and edges $2n \mapsto n,$ $2n+1 \mapsto 3n+2$,
This is an infinite, possibly disconneted graph, call it the 
Collatz graph. The $3x+1$ Conjecture for this
graph can  be formulated as a closed formula in monadic second order logic.
The author considers more generally {\em congruential systems}
which give infinite
graphs with vertices labelled by positive integers 
and edges specified by a finite set of rules $pn+r \mapsto qn+s$,
with $0\le r<p, 0\le s < q$.

In section 1.4 the author  shows that a congruential system without remainders
(all rules have $r= s=0$) can be encoded as a Petri
net with one free place. The reachability problem for a Petri net with
initial position $n$ given is known to be a decidable problem.

A large class of infinite graphs called regular graphs have a decidable
monadic second order theory, see Courcelle [{\em Handbook of Theoretical
Computer Science, Volume B,} Elsevier: New York 1990, pp. 193-242.]
In section 3 the author shows that every regular graph with
bounded vertex degrees can be obtained as a graph of some congruential system.
In section 4 the author shows that the Collatz graph, first viewed as an
undirected graph, and then directed with any choice of orientations of the edges,
is never a regular graph.

\item
Jean Paul Van Bendegem (2005),
{\em The Collatz Conjecture: A Case Study in Mathematical Problem Solving}, 
Logic and Logical Philosophy {\bf 14}, No. 1 (2005), 7--23. (MR 2163301)\\
\newline
\hspace*{.25in}
This philosophical essay concerns the issues of what mathematicans
do beyond proving theorems. The work on the $3x+1$ problem is discussed
from this viewpoint. Such work includes:
computer experiments, heuristic arguments concerning the truth of
the conjecture, metamathemical heuristics concerning the likelihood
of finding a proof, etc.

\item
Stanislav Volkov (2006), 
{\em A probabilistic model for the $5k+1$ problem and related problems,}
Stochastic Processes and Applications {\bf 116} (2006), 662--674. \\
\newline
\hspace*{.25in}
This paper presents a stochastic model for maps like
the $5x+1$ problem, in which most trajectories are
expected to diverge. For the $5x+1$ problem
it  is empirically observed that
the number of values of $n \le x$ that have some
iterate equal to $1$ appears to grow like
$x^{\alpha}$, where $\alpha \approx 0.68$. 
The author develops stochastic models which supply a heuristic
to estimate the value of $\alpha$.

The stochastic model studied is a randomly labelled (rooted)
binary tree model. At each vertex the left branching edge of the tree
gets a label randomly drawn from a (discrete) real distribution X and
the right branching edge gets a label randomly drawn from
a (discrete) real distribution Y. Each vertex is labelled with the
sum of the edge labels from the root; the root gets label $0$. 
The rigorous results of the paper concern such
stochastic models. The author assumes that both $X$ and $Y$ have
positive expected values $\mu_x, \mu_y$, 
but that at least one random variable
assumes some negative values. He also assumes  that the moment generating
functions of both variables are finite for all parameter values.

The author first considers 
for each real $\alpha>0$
the total number of vertices $R_n(\alpha)$ at depth $n$ in
the tree having
label $\le n\alpha$, and sets
$R(\alpha) = \sum_{n=1}^{\infty} R_n(\alpha)$.
He defines a large deviations rate function
$\gamma(\alpha)$ associated to the random variable
$W$ that draws from $-X$ or $-Y$ with equal probability. 
He derives a large deviations criterion (Theorem 1)
which states that if $\gamma(-\alpha) > \log 2$ then $R(\alpha)$
is finite almost surely, while if $\gamma(-\alpha) < \log 2$
then $R(\alpha)$ is infinite almost surely. 
The author next studies the quantity $Q(x)$ counting
the number of vertices with labels smaller than $x$.
This is a refinement of the case $\alpha =0$ above.
He supposes that $\gamma(0) > \log 2$ holds,
and shows (Theorem 2) that this  implies that $Q(x)$ is finite almost surely
for each $x$. He then shows (Theorem 3) that
$$\beta:= \lim_{x \to \infty} \frac{1}{x} \log Q(x)$$
exists almost surely and is given by 
$$
\beta:= \max_{a \in (0, \frac{1}{2}(\mu_x+\mu_y)]} \frac{1}{a}(\log 2 - 
\gamma(-a)).
$$
He constructs a particular  stochastic model of this kind
that approximates the $5x+1$ problem. For this model he shows that
Theorem 3 applies and  computes that $\beta\approx 0.678$.

The author observes that his stochastic
 model has similarities to the branching
random walk stochastic model for
the $3x+1$ problem studied in Lagarias and Weiss (1992),
whose analysis also used the theory of large deviations. 

{\em Note.} The  exponential
branching of the $5x+1$ tree above $1$ allows one to prove that the
number of such $n \le x$ that have some $5x+1$ iterate
equal to $1$ is at least $x^{\beta}$ for
some small positive $\beta$.

\item
Wang, Nian-liang (2002),
{\em The way to prove Collatz problem} (Chinese), 
Journal of Shangluo University (2002), No. 2. \\
\newline
\hspace*{.25in}
[I have not seen this paper.]

\item
Wang, Xing-Yuan; Wang, Qiao Long; Fen, Yue Ping; and Xu, Zhi Wen (2003), 
{\em The distribution of the fixed points on the real axis of a generalized
$3x+1$ function and some related fractal images} (Chinese),
Journal of Image and Graphics [Zhongguo tu xiang tu xing xue bao] {\bf 8}(?) (2003), No. 1.\\
\newline
\hspace*{.25in}
[I have not seen this paper.]

\item
Xing-yuan Wang  and  Xue-jing Yu (2007),
{\em Visualizing generalized $3x+1$ function dynamics based on fractal,}
Applied  Mathematics and Computation 
{\bf 188} (2007), no. 1, 234--243. 
(MR2327110).\\
\newline
\hspace*{.25in}
This paper studies two complex-valued generalizations of
 the Collatz function. It replaces the function
 $mod_2(x)$ defined on integers $x$ by the function $(\sin \frac{pi x}{2})^2$,
 defined for complex $x$. These are then substituted in 
 the definitions 
 $$
 C(x) = \frac{x}{2}( 1 - mod_2(x)) + (3x+1) mod_2(x)
 $$
 and 
 $$
 T(x) = \frac{3^{mod_2(x)} + mod_2(x)}{2^{1-mod_2(x)}}.
 $$
 Both these functions agree with the Collatz function $C(n)$ on
 the positive integers.  The first function simplifies to 
 $$
 C(x) = \frac{1}{4}(7x+2 - (5x+2) \cos \pi x ) .
 $$
 The paper  studies 
iteration of the entire functions $C(x)$ and $T(x)$,  from the viewpoint
of escape time, stopping time and total stopping time. It presents graphics
illustrating the results of the algorithms. The total stopping time plots  exhibit 
vaguely Mandlebrot-like sets of various sizes located around the positive integers.

\item
Wang, Xing-yuan and Yu, Xue-jing (2008),
{\em Dynamics of the generalized $3x+1$ function determined by
its fractal images,}
Progress Natural Science (English Edition) {\bf 18} (2008), 217--223. (MR 2009i:37020).\\
\newline
\hspace*{.25in}
From the abstract:  ``Two different complex maps are obtained
by generalizing the $3x+1$ function to the complex plane, and fractal
images for these two complex maps are constructed by using escape-time,
stopping time and total-stopping-time arithmetic."

\item
G\"{u}nther J. Wirsching (2000),
{\em \"{U}ber das $3n + 1$ Problem,}
Elem. Math. {\bf 55} (2000), 142--155. 
(MR 2002h:11022). \\
\newline
\hspace*{.25in}
This is a  survey paper, which discusses
the origin of the $3n+1$ problem and 
results on the dynamics of the $3x+ 1$ function.

\item
G\"{u}nther J. Wirsching (2001),
{\em A functional differential equation and $3n + 1$ dynamics,}
in: {\em Topics in Functional Differential and Functional Difference
Equations (Lisbon 1999),} (T. Faria, E. Frietas, Eds.), Fields Institute
Communications No. 29, Amer. Math. Soc. 2001, pp. 369--378. 
(MR 2002b:11035). \\
\newline
\hspace*{.25in}
This paper explains how a functional differential equation
arises in trying to understand $3n + 1$ dynamics. as given
in Wirsching (1998a). It analyzes some properties of its
solutions.

\item
G\"{u}nther  J. Wirsching (2003)
{\em On the problem of positive predecessor density in $3N+1$ dynamics,}
Disc. Cont. Dynam. Syst. {\bf 9} (2003), no. 3, 771--787.
(MR 2004f:39028). \\
\newline
\hspace*{.25in}
This paper discusses an approach to prove positive predecessor
density, which formulates three conjectures which, if proved,
would establish the result. This approach presents in more detail
aspects of
the approach taken in the author's Springer Lecture Notes volume,
Wirsching (1998a).

\item
Wu, Jia Bang and  Huang, Guo Lin (2001),
{\em On the traditional iteration and the elongate iteration of the $3N+1$ conjecture} (Chinese),
Journal of South-Central University for Nationalities. Natural Sciences Ed.
[Zhong nan min zu  da xue xue bao. Zi ran ke xue   ban] 
{\bf 21} (2001), No. 4 \\
\newline
\hspace*{.25in}
[I have not seen this paper.]

\item
Wu, Jia Bang  and Hao, Shen Wang (2003),
{\em On the equality of the adequate stopping time and the coefficient stopping
time of $n$ in the $3N+1$ conjecture} (Chinese, English summary),
J. Huazhong Univ. Sci. Technol. Nat. Sci.  [Hua zhong gong xue yuan]  {\bf 31} (2003), no. 5, 114--116. (MR 2000420)\\
\newline
\hspace*{.25in}
English summary: "In the paper the equality of adequate stopping time $t_a(n)$ and the
coefficient stopping time $t_c(n)$ in the $3N+1$ conjecture was discussed. It was
proved that $t_a(n) = t_c(n)$ on condition that $d= \sum_{i=1}^{k-1} x_i(n)$ is not
very large. Therefore it was conjectured that when the bound condition for $d$
was unnecessary, or was automatically satisfied, $t_c(n)= t_a(n)$ in all cases."

{\em Note:} The notion of coefficient stopping time was introduced in Terras (1976)
for the $3x+1$ function;
here $t_c(n)$ is its analogue for the Collatz function.
The adequate stopping time $t_a(n)$ is the analogue of the stopping time for
the Collatz function. Terras (1976) made a conjecture equivalent to
asserting  the equality $t_a(n)= t_c(n)$ always holds.

\item
Wu, Jia Bang and Huang, Guo Lin  (2001a),
{\em Family of consecutive integer pairs of the same height in the Collatz conjecture}
(Chinese, English summary), Mathematica  Applicata (Wuhan) [Ying yung shu hs\"{u}eh] {\bf 14} (2001), suppl. 21--25.
(MR 1885838, Zbl 1002.11023)\\
\newline
\hspace*{.25in}
English summary: "In the Collatz conjecture, certain runs of consecutive integers have 
the same height. In particularly in this paper, pairs of successive integers of the same
height are investigated. It is found that families of consecutive integer pairs of the same
height occur infinitely often and in different patterns."

{\em Note.} The authors show coalescence of iteration of the $3x+1$ function for consecutive
pairs in arithmetic progessions  of form $\{(2^k m + r, 2^k m + r+1): m\ge 0\}$,
for example $(32m+5, 32m+6), (64m+45, 64m+46), (128m+29, 128m+30)$.
They exhibit infinitely many such arithmetic progressions, specified by the pattern of
residues (modulo 2) of the successive iterates of the $3x+1$
function, e.g. for $(128m+29, 128m+30)$
the  patterns are $(1001101), (0111100)$,  respectively.

\item
Wu, Jia Bang  and Huang, Guo Lin (2001b),
{\em Elongate iteration for the $3N+1$ Conjecture}
(Chinese, English summary), J. Huazhong Univ. Sci. Tech.  
[Hua zhong gong xue yuan]  {\bf 29} (2001), no.2 , 112--114.
(MR 1887558)\\
\newline
\hspace*{.25in}
English summary: "The concept of Elongate Iteration for the $3N+1$ conjecture is
presented and discussed. Some results of Elongate Iteration are given: {\bf a.}
The correspondence of numbers with parity vectors.
{\bf b.} Invariability of $l$-tuple. {\bf c.} Proof of term formula of $n$. {\bf d.} Some
equivalence propositions for the $3N+1$ conjecture. {\bf e.} A property of
coefficient stopping time $t_c(n)$, etc."

\item
Wu, Wen Quan (2003),
{\em An equivalent proposition of Collatz conjecture} (Chinese),
Journal of Aba Teacher's College 
[Aba shi fan gao deng zhuan ke xue xiao xue bao] (2003), No. 3.\\
\newline
\hspace*{.25in}
[I have not seen this paper.]

\item
Xia, Ye (2003),
{\em Some thoughts about $3X+1$ Conjecture}  (Chinese),
[Zhongxue Shuxue Za Zhi (Chu zhong ban)] (2003), No. 6. \\
\newline
\hspace*{.25in}
[I have not seen this paper.]

\item
Xie, Jiu Guo  (2006),
{\em The interesting $3x+1$ problem (Chinese), }
J. of Hunan Univ. Science Technology 
[Hunan ki ji da xue xue bao] (2006), No. 11 \\
\newline
\hspace*{.25in}
[I have not seen this paper.]

\item
Michinori Yamashita (2002),
{\em  $3x+1$ problem from the $(e,k)$-perspective } (Japanese),
PC Literacy [Pasocon Literacy] (Personal Computer Users Application Technology
Association]  {\bf 27} (2002), No. 10,  22--27.\\
\newline
\hspace*{.25in}
[I have not seen this paper.]

\item
Michinori Yamashita (2006),
{\em Note on the $3x \pm 1$ Problem } (Japanese),
PC Literacy [Pasocon Literacy] (Personal Computer Users Application Technology
Association]  (2006).\\
\newline
\hspace*{.25in}
This paper studies the $3x+1$ function $f(x)=T(x)$ and the $3x-1$ function
$g(x)= -T(-x)$. It lets $(e,k)=2^ek -1,~[e,k]= 2^e k +1$, where $k$ is an odd integer.
It proves a number of identities for these functions, such as
$3 (e,k) =(1, (e-1, 3k) ),$  given a table of what happens when one multiplies
by a small number $l$.  It studies patterns of successive iterations.


\item
Roger E.  Zarnowski (2001),
{\em Generalized inverses and the total stopping time of
Collatz sequences,}
 Linear and Multilinear Algebra {\bf 49} (2001), 115--130. 
(MR 2003b:15011). \\
\newline
\hspace*{.25in}
The $3x+1$ iteration is formulated in terms of a
denumberable Markov chain with transition matrix $P$.
The $3x+1$ Conjecture is reformulated in terms of the
limiting behavior of $P^{k}$. The group inverse
$A^{\sharp}$ to an $n \times n$ matrix $A$ is 
defined by the properties $A A^{\sharp} =  A^{\sharp} A$,
$A A^{\sharp} A = A$ and $A^{\sharp} A A^{\sharp} = 
A^{\sharp} $, and is unique when it exists.
Now set  $A = I - P$, an infinite matrix.
Assuming there are
no nontrivial cycles, the group inverse
$A^{\sharp}$ exists,
and satisfies  $\lim_{k \to \infty} P^k = I - A A^{\sharp}.$
An explicit formula is given for $A^{\sharp}$.

\item
Roger E. Zarnowski (2008),
{\em The congruence structure of the $3x+1$ map,}
Fibonacci Quarterly {\bf 46/47} (2008/09), no. 2, 115--125.
(MR 2010c:11037)\\
\newline
\hspace*{.25in}
This paper studies the structure of forward iteration of
the $3x+1$ map taking congruence classes $(\bmod~2^n)$ to
congruence classes $(\bmod ~3^k)$ where $k$ depends on
the parity sequence of the iterates. He describes relations of
image classes (which may overlap), using a concept termed
``congruence class triangle."

This paper observes that the set $A$ of nonnegative integers 
relatively prime to $3$ is an invariant set under forward iteration
of the $3x+1$ map. It then shows that under iteration 
the set of positive integers in any fixed congruence class
$(\bmod ~2^a 3^b)$ has forward orbits visiting every integer in $A$
infinitely often. Thus to prove the $3x+1$ conjecture it suffices to
verify it on integers in any one of these congruence classes. 

{\em Note.} There are quite a few prior results known on
sufficiency of $3x+1$ conjecture on  a ``thin" congruence
implying it is true in general, e.g.  Korec and Znam (1987).

\item
Zhou, Yun Cai and Zhou, Bao Lan (2007),
{\em On the Collatz Conjecture} (Chinese),
J. Yangtze University Natural Science.
[Changjiang daxue xuebao. Zi ke ban] (2007), No. 1, 24--26.\\
\newline
\hspace*{.25in}
English summary: ``The Collatz wave set (CWS) and the 3-1 mapping trim applied
on it are defined, and the properties of CWS under trim mapping are discussed.
Collatz's conjecture, namely, the $3x+1$ problem, is also discussed from the
viewpoint of CWS and 3-1 mapping, by which a new possibility is provided for
proving Collatz's conjecture."

\end{enumerate}

\noindent
{\bf Acknowledgements.}
I thank  M. Chamberland, 
J. Goodwin, C. Hewish,
S. Kohl, Wang Liang,  C. Reiter,  J. L. Simons, T. Urata, and S. Volkov and 
G. J. Wirsching for supplying references.  I thank
W.-C. Winnie Li, Xinyun Sun and Yang Wang for help with Chinese references.
I thank  T. Ohsawa for help with Japanese references.  \\


\end{document}